\documentclass[final,3p,times,12pt]{elsarticle}
\usepackage{amsfonts}
\usepackage{bbm}
\usepackage{txfonts}
\usepackage{}
\usepackage{amssymb}
\usepackage{amsmath}
\usepackage{graphicx}
\usepackage{subfigure}

\usepackage{array}
\usepackage{mathrsfs}
\usepackage{mathptmx}
\usepackage{color}

\newtheorem{Theo}{Theorem}

\newtheorem{Rem}{Remark}
\newtheorem{Lem}{Lemma}

\numberwithin{equation}{section}
\numberwithin{Lem}{section}
\numberwithin{Defi}{section}
\numberwithin{Theo}{section}
\numberwithin{Pro}{section}
\numberwithin{Rem}{section}
\numberwithin{Coro}{section}
\numberwithin{Fig}{section}

\begin{document}
\def\b{\beta}
\def\a{\alpha}
\def\o{\"{o}}
\def\t{\tau}
\begin{frontmatter}

\title{An implicit midpoint difference scheme for the fractional Ginzburg-Landau equation 
}
\author{Pengde Wang, 
 Chengming Huang$^{*}$ }

\cortext[cor1]{Corresponding author.\\
\emph{Email addresses}: \texttt{pengde\_wang@yeah.net} (P. Wang),
\texttt{chengming\_huang@hotmail.com} (C. Huang)}
\address{ School of Mathematics and Statistics, Huazhong University of Science and
Technology, Wuhan 430074, China\\}

\date{}
\begin{abstract}
\par
This paper proposes and analyzes an efficient difference scheme for
the nonlinear complex Ginzburg-Landau equation involving fractional
Laplacian. The scheme is based on the implicit midpoint rule for the
temporal discretization and a weighted and shifted Gr\"{u}nwald
difference operator for the spatial fractional Laplacian. By virtue
of a careful analysis of the difference operator, some useful
inequalities with respect to suitable fractional Sobolev norms are
established. Then the numerical solution is shown to be bounded, and
convergent in the $l^2_h$ norm with the optimal order $O(\t^2+h^2)$
with time step $\t$ and mesh size $h$. The a priori bound as well as
the convergence order hold unconditionally, in the sense that no
restriction on the time step $\t$ in terms of the mesh size $h$
needs to be assumed. Numerical tests are performed to validate the
theoretical results and effectiveness of the scheme.
\end{abstract}
\begin{keyword}
Fractional Ginzburg-Landau  equation; Fractional Laplacian; Riesz fractional derivative; Weighted and shifted Gr\"{u}nwald difference; Convergence
\end{keyword}

\end{frontmatter}


\section{Introduction}
%
The classical complex Ginzburg-Landau equation (GLE) is one of the
most-studied nonlinear equations in the physics community, which
describes a vast variety of phenomena from nonlinear waves to
second-order phase transitions, from superconductivity,
superfluidity, and Bose-Einstein condensation to liquid crystals and
strings in field theory \cite{AransonIS:2002:RMP}. The fractional
generalization of the GLE was suggested in
\cite{TarasovV:2005:PA,TarasovV:2006:C} from the variational
Euler-Lagrange equation for fractal media. Since then, the
fractional Ginzburg-Landau equation (FGLE) has been exploited to
describe various physical phenomena, such as the dynamical processes
in continuums with fractal dispersion and the media with fractal
mass dimension \cite{TarasovV:2005:PA}, a fairly general class of
critical phenomena when the organization of the system near the
phase transition point is influenced by a competing nonlocal
ordering \cite{MilovanovA:2005:PLA} and a network of diffusively
Hindmarsh-Rose neurons with long-range synaptic coupling
\cite{MvogoA:2014:Arxiv}.

In this paper, we consider the following FGLE with the fractional Laplacian $(1<\a\leqslant 2)$
\begin{equation}\label{eq:FSE1}
  u_t+(\upsilon+i\eta)(-\Delta)^{\frac{\a}{2}}u+(\kappa+i\zeta){|u|}^2u-\gamma u=0,\quad x\in\mathbb{R},\ t\in(0,T],
\end{equation}
subject to the initial condition
\begin{equation}\label{eq:IC1}
  u(x,0)=u_0(x), \quad x\in\mathbb{R},
\end{equation}
where $i=\sqrt{-1}$, $u(x,t)$ is a complex-valued function of time $t$ and space $x$, $\upsilon>0, \kappa>0$, $\eta, \zeta, \gamma$ are given real constants and $u_0(x)$ is a given smooth function. The
fractional Laplacian can be regarded as a pseudo-differential operator
with the symbol $-|\xi|^{\a}$:
\begin{equation*}
  -(-\Delta)^{\frac{\a}{2}}u(x,t)=-\mathcal{F}^{-1}(|\xi|^{\a}\hat{u}(\xi,t)),
\end{equation*}
where $\mathcal{F}$ denotes the Fourier transform.
It is indeed equivalent to the
following Riesz fractional derivative \cite{ZhuangP:2009:SIAMNUA,YangQ:2010:AMM}, i.e.,
\begin{equation}\label{eq:FDef}
   -(-\Delta)^{\frac{\a}{2}}u(x,t)=\frac{{\partial}^{\a}}{\partial{|x|}^{\a}}u(x,t)
   :=-\frac{1}{2\cos{\frac{\a\pi}{2}}}\big[{}_{-\infty}D^{\a}_{x}u(x,t)+{}_{x}D^{\a}_{+\infty}u(x,t)\big],
\end{equation}
where ${}_{-\infty}D^{\a}_{x}u(x,t)$ denotes the left Riemann-Liouville fractional derivative \cite{PodlubnyI:1999:AP}
\[
{}_{-\infty}D^{\a}_{x}u(x,t)=\frac{1}{\Gamma(2-\a)}\frac{d^2}{dx^2}\int^x_{-\infty}\frac{u(\xi,t)}{(x-\xi)^{\a-1}}d\xi,
\]
and ${}_{x}D^{\a}_{+\infty}u(x,t)$ the right Riemann-Liouville fractional derivative
\[
{}_{x}D^{\a}_{+\infty}u(x,t)=\frac{1}{\Gamma(2-\a)}\frac{d^2}{dx^2}\int^{+\infty}_x\frac{u(\xi,t)}{(\xi-x)^{\a-1}}d\xi.
\]
Obviously, when $\a=2$, this operator reduces to the classical Laplacian and the equation reduces to the classical
cubic nonlinear complex GLE. 

From the mathematical point of view, unlike the case of classical GLE, the
dissipative mechanism of the FGLE is not characterized by the classical Laplacian but by the
fractional power of the Laplacian, which has raised some essential difficulties in theoretical analysis
and recently drawn quite a lot of interest from various authors. For example, Tarasov \cite{TarasovV:2006:JPA} derived and analyzed
the psi-series solution. 
Pu and Guo \cite{PuX:2013:AA} investigated the global well-posedness, long-time dynamics and global attractors for the nonlinear FGLE. Guo ang Huo \cite{GuoB:2013:FCAA} studied the inviscid limit behavior of the FGLE to the fractional Schr\o dinger equation (FSE). Lu et al. \cite{LuS:2013:IJBC} analyzed the well-posedness and asymptotic behaviors in two dimensions.
Millot and Sire \cite{MillotV:2015:ARMA} considered the asymptotic analysis of the FGLE in bounded domain and showed that solutions with uniformly
bounded energy converge weakly to sphere valued 1/2-harmonic maps.

From the numerical point of view, however, there is very little attention to the numerical solution of the FGLE.
 In order to simulate the propagation of the localized
impulses in diffuse neural networks, Mvogo et al.
\cite{MvogoA:2014:Arxiv} proposed a semi-implicit Riesz fractional
finite difference scheme, which is first order in time and second
order in space. Numerical simulations show that the scheme is
feasible and efficient. To the best of our knowledge, it seems that
this is the only work. Nevertheless, that paper focused on the establishment and simulation of the model equation.
The stability and convergence of the discretization scheme were not
discussed.

In this paper, we propose another difference scheme for the FGLE,
which treats the time derivative by the implicit midpoint rule and
the space derivative by the second-order accurate weighted and
shifted Gr\"{u}nwald difference (WSGD) method \cite{HaoZ:2015:JCP}.
This scheme is second-order in both time and space. Our focus is on
a rigorous theoretical analysis for the scheme. We will prove that
the scheme is unconditionally convergent with optimal order, in the
sense that no added restriction on the temporal step size in terms
of the spatial discretization parameter needs to be assumed. We
mention that there are a number of convergence results for the
classical GLE in the literature (see, e.g.,
\cite{LordGJ:1997:SIAMNU,SunZ:1998:JCAM,MatsuoT:2001:JCP,XuQ:2011:NMPDE,WangT:2011:NMPDE,ZhangL:2010:NFAO,ZhangY:2013:NMPDE,HaoZ:2015:NMPDE}
and reference therein), which are normally derived based on the
uniform boundedness of numerical solutions
\cite{LordGJ:1997:SIAMNU,SunZ:1998:JCAM}. In our case, however, it
seems difficult for us to follow this approach because of the
nonlocal property of the fractional Laplacian. In order to overcome
this obstacle, we make a detailed study of the fractional
approximation operator. The discrete fractional Gagliardo-Nirenberg
inequality and an equivalence relation between an energy norm and the
fractional Sobolev semi-norm are established.

The remainder of this paper is arranged as follows. In Section 2, we introduce
the WSGD operator to discretize the
involved fractional derivative, and give some 
technical
lemmas. In Section 3, we 
establish our fully discrete scheme. Section 4 is devoted to the
rigorous theoretical analysis, including unique solvability and
especially, the boundedness and convergence. In Section 5, we carry
out some numerical experiments to confirm our theoretical results
and show the efficiency of the proposed scheme. Finally, we  draw
some conclusions in Section 6.
%
\section{Preliminaries}

\subsection{Spatial discretization}
Up to now, a broad range of difference methods have been constructed
to approximate the Riemann-Liouville fractional derivative and Riesz
fractional derivative, such as the Gr\"{u}nwald-based scheme
\cite{MeerschaertM:2004:JCAM,TadjeranC:2006:JCP,TianW:2015:MC,ZhouH:2013:JSC,HaoZ:2015:JCP,ChenM:2014:SIAMNUA,ChenM:2014:CICP,BaeumerB:2015:TAMS,ZhaoL:2015:NMPDE}
and the fractional centered difference based scheme
\cite{Ortigueira:2006:IJMMS,CelikC:2012:JCP,DingH:2015:JCP,ZhaoX:2014:SIAMSC}.
Each of them has its own advantages.

In this paper, we use the WSGD method, proposed by Hao, Sun and Cao
\cite{HaoZ:2015:JCP}, to approximate the left and right
Riemann-Liouville space fractional derivatives, respectively. The
essential idea of this approximation is using the weighted average
to vanish the low order leading terms in asymptotic expansions for
the truncation errors of the shifted Gr\"{u}nwald formulae. Another
class of WSGD approximations can be found in Tian, Zhou and Deng
\cite{TianW:2015:MC}. Initially, the shifted Gr\"{u}nwald formulae
were constructed by Meerschaert and Tadjeran
\cite{MeerschaertM:2004:JCAM} and defined by
\begin{equation}
\begin{split}
&{}_{L}\mathcal{A}^{\a}_{h,p}u(x):=\frac{1}{h^{\a}}\sum^{\infty}_{l=0}g^{(\a)}_lu(x-(l-p)h)={}_{-\infty}D^{\a}_{x}u(x)+O(h),\\
&{}_{R}\mathcal{A}^{\a}_{h,r}u(x):=\frac{1}{h^{\a}}\sum^{\infty}_{l=0}g^{(\a)}_lu(x+(l-r)h)={}_{x}D^{\a}_{+\infty}u(x)+O(h),
\end{split}
\end{equation}
where $p, r$ are integers and $g^{(\a)}_l=(-1)^l\binom{\a}{l}$ are the coefficients of the power series of the function $(1-z)^{\a}$, i.e.,
\begin{equation}\label{eq:GL}
(1-z)^{\a}=\sum^{\infty}_{l=0}(-1)^l\binom{\a}{l}z^l=\sum^{\infty}_{l=0}g^{(\a)}_lz^l,
\end{equation}
for $|z|< 1$, and they can be evaluated recursively
\[
g^{(\a)}_0=1, \quad g^{(\a)}_l=(1-\frac{\a+1}{l})g^{(\a)}_{l-1}, \quad l=1,2,\ldots.
\]
By weighting the Gr\"{u}nwald
approximation formulae with different shifts, Hao, Sun and Cao
\cite{HaoZ:2015:JCP} 
propose the following WSGD operator:
\begin{equation}
\begin{split}
&{}_{L}\mathcal{D}^{\a}_{h}u(x)=\lambda_1{}_{L}\mathcal{A}^{\a}_{h,1}u(x)+\lambda_0{}_{L}\mathcal{A}^{\a}_{h,0}u(x)+\lambda_{-1}{}_{L}\mathcal{A}^{\a}_{h,-1}u(x),\\
&{}_{R}\mathcal{D}^{\a}_{h}u(x)=\lambda_1{}_{R}\mathcal{A}^{\a}_{h,1}u(x)+\lambda_0{}_{R}\mathcal{A}^{\a}_{h,0}u(x)+\lambda_{-1}{}_{R}\mathcal{A}^{\a}_{h,-1}u(x),
\end{split}
\end{equation}
where
\begin{equation}\label{eq:c1}
\lambda_1=\frac{\a^2+3\a+2}{12},\quad
\lambda_0=\frac{4-\a^2}{6},\quad \lambda_{-1}=\frac{\a^2-3\a+2}{12}.
\end{equation}

%

\begin{Lem}\label{lem:WSGD} (See \cite{HaoZ:2015:JCP}.)
Suppose that $u\in L^1(\mathbb{R})$ and 
\[
u\in\mathscr{L}^{2+\a}(\mathbb{R}):=\Big\{u\Big|\int^{+\infty}_{-\infty}\big(1+|\xi|\big)^{2+\a}|\widehat{u}(\xi)|\text{d}\xi<\infty\Big\},
\]
then for a fixed h, we have
\begin{equation}
\begin{split}
&{}_{L}\mathcal{D}^{\a}_{h}u(x)={}_{-\infty}D^{\a}_{x}u(x)+O(h^2),\\
&{}_{R}\mathcal{D}^{\a}_{h}u(x)={}_{x}D^{\a}_{+\infty}u(x)+O(h^2).
\end{split}
\end{equation}
\end{Lem}

Rearranging the WSGD operator gives
\begin{equation}\label{eq:WSGDR1}
\begin{split}
&{}_{L}\mathcal{D}^{\a}_{h}u(x):=\frac{1}{h^{\a}}\sum^{\infty}_{l=0}w^{(\a)}_lu(x-(l-1)h)={}_{-\infty}D^{\a}_{x}u(x)+O(h^2),\\
&{}_{R}\mathcal{D}^{\a}_{h}u(x):=\frac{1}{h^{\a}}\sum^{\infty}_{l=0}w^{(\a)}_lu(x+(l-1)h)={}_{x}D^{\a}_{+\infty}u(x)+O(h^2),
\end{split}
\end{equation}
where
\begin{equation}\label{eq:cpq}
\left\{\begin{array}{ll}
   w^{(\a)}_0=\lambda_1g^{(\a)}_0,\quad w^{(\a)}_1=\lambda_1g^{(\a)}_1+\lambda_0g^{(\a)}_0,\\
   w^{(\a)}_l=\lambda_1g^{(\a)}_l+\lambda_0g^{(\a)}_{l-1}+\lambda_{-1}g^{(\a)}_{l-2},\quad l\geqslant2.
     \end{array}\right.
\end{equation}
In addition, the coefficients have the following properties:
\begin{equation}\label{eq:cpq1}
\left\{\begin{array}{ll}
   w^{(\a)}_0\geqslant0,\quad w^{(\a)}_1\leqslant0,\quad w^{(\a)}_l\geqslant0,\quad l\geqslant3,\\
   \sum^{+\infty}_{l=0}w^{(\a)}_l=0,\quad \sum^{M}_{l=0}w^{(\a)}_l\leqslant0,\quad M\geqslant1,\\
   w^{(\a)}_0+w^{(\a)}_2\geqslant0.
     \end{array}\right.
\end{equation}
\begin{Rem}
For $1<\a<2$, it is easy to verify that the inequalities in \eqref{eq:cpq1} are strictly true, i.e., the sign $``\leqslant"$ and $``\geqslant"$ can be substituted by $``<"$ and $``>"$.
\end{Rem}

Using the relation \eqref{eq:FDef} and \eqref{eq:WSGDR1}, the WSGD
approximation for the fractional Laplacian can be given by
\begin{equation}\label{eq:wsgdr1}
\begin{split}
\Delta^{\a}_{h}u(x)&:=\frac{1}{2\cos{\frac{\a\pi}{2}}}\Big({}_{L}\mathcal{D}^{\a}_{h}u(x)+{}_{R}\mathcal{D}^{\a}_{h}u(x)\Big)\\
&=\frac{1}{h^{\a}}\frac{1}{2\cos{\frac{\a\pi}{2}}}\Big(\sum^{\infty}_{l=0}w^{(\a)}_lu(x-(l-1)h)+\sum^{\infty}_{l=0}w^{(\a)}_lu(x+(l-1)h)\Big)\\
&=(-\Delta)^{\frac{\a}{2}}u(x)+O(h^2).
\end{split}
\end{equation}
\subsection{Fractional Sobolev norm}
Now we introduce some fractional Sobolev norms and relevant lemmas.

Let $\mathbb{Z}$ denote the set of all integers
and $h\mathbb{Z}$ denote the infinite grid with grid points $x_j=jh$
for $j\in\mathbb{Z}$. For any grid functions $u=\{u_j\}$,
$v=\{v_j\}$ on $h\mathbb{Z}$, we define the discrete inner product
and the associated $l^2_h$ norm as
\[
(u,v)_h=h\sum_{j\in\mathbb{Z}}u_j\overline{v}_j,\quad
\|u\|^2_h=(u,u)_h.
\]
We also define the discrete $l^p_h$ norm as
\[
\|u\|^p_{l^p_h}=h\sum_{j\in\mathbb{Z}}|u_j|^p, \quad 1\leqslant p<+\infty,
\]
and the discrete $l^{\infty}_h$ norm as
\[
\|u\|_{l^\infty_h}=\sup_{j\in\mathbb{Z}}|u_j|.
\]

Define the space $l^2_h:=\{u \ | \ u=\{u_j\},
\|u\|_h<+\infty\}$. For any $u \in l^2_h$, the semi-discrete Fourier
transform of $u$ is the function $\widehat{u}\in L^2[-\pi/h,\pi/h]$
defined by
\[
\widehat{u}(k):=\frac{1}{\sqrt{2\pi}}h\sum_{j\in\mathbb{Z}}u_je^{-ikx_j};
\]
see \cite{TrefethenLN:2000:SMS,StrikwerdaJC:2004:FDSPDE}. Moreover, we have the inversion formula
\[
u_j=\frac{1}{\sqrt{2\pi}}\int^{+\pi/h}_{-\pi/h}\widehat{u}(k)e^{ikx_j}\text{d}k,
\]
and Parseval's theorem gives
\begin{equation}\label{eq:PI}
(u,v)_h=\int^{+\pi/h}_{-\pi/h}\widehat{u}(k)\overline{\widehat{v}(k)}\text{d}k.
\end{equation}
Given a constant $\sigma\in [0, 1]$, we define the fractional
Sobolev semi-norm $|\cdot|_{H^\sigma_h}$  and norm
$\|\cdot\|_{H^\sigma_h}$ as
\begin{equation}
|u|_{H^\sigma_h}^2=\int^{+\pi/h}_{-\pi/h}|k|^{2\sigma}|\widehat{u}(k)|^2\text{d}k,\quad
\|u\|_{H^\sigma_h}^2=\int^{+\pi/h}_{-\pi/h}(1+|k|^{2\sigma})|\widehat{u}(k)|^2\text{d}k.
\end{equation}
From \eqref{eq:PI}, it is clear that $\|u\|_{H^\sigma_h}^2=\|u\|^2_h+|u|_{H^\sigma_h}^2$ and $|u|_{H^0_h}^2=\|u\|^2_h$.
Then we introduce the following lemmas.

\begin{Lem} 
For $0\leqslant \sigma_0\leqslant \sigma\leqslant 1$, there exists a
constant $C\in[1, \sqrt{2}]$ such that
\begin{equation}\label{eq:HHi}
\|u\|_{H^{\sigma_0}_h}\leqslant
C\|u\|_{H^{\sigma}_h}^{\frac{\sigma_0}{\sigma}}\|u\|_h^{1-{\frac{\sigma_0}{\sigma}}}.
\end{equation}
\end{Lem}
\textbf{Proof.} From the definition of $\|u\|_{H^{\sigma_0}_h}$ and
H\o lder's inequality, we have
\begin{equation}
\begin{split}
\|u\|^2_{H^{\sigma_0}_h}&=\int^{+\pi/h}_{-\pi/h}(1+|k|^{2\sigma_0})|\widehat{u}(k)|^2\text{d}k\\
&=\int^{+\pi/h}_{-\pi/h}\big((1+|k|^{2\sigma})|\widehat{u}(k)|^2\big)^{\frac{\sigma_0}{\sigma}}
\big(|\widehat{u}(k)|^2\big)^{1-\frac{\sigma_0}{\sigma}}
\Big(\frac{1+|k|^{2\sigma_0}}{(1+|k|^{2\sigma})^{\frac{\sigma_0}{\sigma}}}\Big)\text{d}k\\
&\leqslant C^2\Big(\int^{+\pi/h}_{-\pi/h}(1+|k|^{2\sigma})|\widehat{u}(k)|^{2}\text{d}k\Big)^{\frac{\sigma_0}{\sigma}}
\Big(\int^{+\pi/h}_{-\pi/h}|\widehat{u}(k)|^{2}\text{d}k\Big)^{1-\frac{\sigma_0}{\sigma}}\\
&=C^2\Big(\|u\|_{H^{\sigma}_h}^{\frac{\sigma_0}{\sigma}}\|u\|_h^{1-{\frac{\sigma_0}{\sigma}}}\big)^2,
\end{split}
\end{equation}
where we have used the inequality
$\frac{1}{2}(1+a^\mu)\leqslant(1+a)^\mu\leqslant(1+a^\mu)$ for $a>0,
0\leqslant\mu\leqslant1$ to derive the third line of above
inequality. Thus the proof is complete. $\Box$
%
%
%
\begin{Lem}\label{lem:dgni}
For any $\frac{p-2}{2p}<\sigma_0\leqslant 1$, there exists a
constant $C_{\sigma_0}=C(\sigma_0)>0$ independent of $h>0$, such that
\begin{equation}\label{eq:dgni}
\|u\|_{l^p_h}\leqslant C_{\sigma_0}\|u\|^{\frac{\sigma_0}{\sigma}}_{H^{\sigma}_h}
\|u\|_h^{1-\frac{\sigma_0}{\sigma}},\quad 2\leqslant p\leqslant +\infty,
\end{equation}
for every $\sigma_0\leqslant \sigma \leqslant 1$.
\end{Lem}
\textbf{Proof.}
Using the Hausdorff-Young inequality (see Appendix A), for $1\leqslant q\leqslant 2$ such that $\frac{1}{p}+\frac{1}{q}=1$, we have
\begin{equation}
\begin{split}
\Big(h\sum_{j\in\mathbb{Z}}|u_j|^p\Big)^{\frac{1}{p}}
&\leqslant C\Big(\int^{+\pi/h}_{-\pi/h}|\widehat{u}(k)|^{q}\text{d}k\Big)^{\frac{1}{q}}\\
&=C\Big(\int^{+\pi/h}_{-\pi/h}\frac{1}{\big(1+|k|^{2\sigma_0}\big)^{\frac{q}{2}}}\big(1+|k|^{2\sigma_0}\big)^{\frac{q}{2}}|
\widehat{u}(k)|^{q}\text{d}k\Big)^{\frac{1}{q}}.
\end{split}
\end{equation}
From the H\o lder's inequality, it follows that
\begin{equation}
\begin{split}
\Big(h\sum_{j\in\mathbb{Z}}|u_j|^p\Big)^{\frac{1}{p}}
&\leqslant C\Big(\int^{+\pi/h}_{-\pi/h}\big(1+|k|^{2\sigma_0}\big)|
\widehat{u}(k)|^2\text{d}k\Big)^{\frac{1}{2}}
\Big(\int^{+\pi/h}_{-\pi/h}\frac{1}{\big(1+|k|^{2\sigma_0}\big)^{\frac{q}{2-q}}}\text{d}k\Big)^{\frac{2-q}{2q}}\\
&\leqslant C\|u\|_{H^{\sigma_0}_h}\Big(\int^{+\infty}_{-\infty}\frac{1}{\big(1+|z|^{2\sigma_0}\big)^{\frac{q}{2-q}}}
\text{d}z\Big)^{\frac{2-q}{2q}}.
\end{split}
\end{equation}
Then for $\frac{p-2}{2p}<\sigma_0\leqslant 1$, we obtain
\begin{equation*}
\|u\|_{l^p_h}\leqslant \tilde{C}_{\sigma_0}\|u\|_{H^{\sigma_0}_h},
\end{equation*}
where $\tilde{C}_{\sigma_0}=\tilde{C}(\sigma_0)>0$ is independent of
$h$. Combining above inequality with
\eqref{eq:HHi} gives \eqref{eq:dgni} and thus completes the proof.
$\Box$
\begin{Rem}
Lemma \ref{lem:dgni} is an extension of Lemma 3.2 in
\cite{KirkpatrickK:2013:CMP}, where the special case with $p=4$ is
considered.
\end{Rem}
%
%
\begin{Lem}\label{Lem:L1}
For $1<\a\leqslant 2$, let $h(\a,\omega)$ be the function defined by
\begin{equation}\label{eq:hak}
h(\a,\omega)=\lambda_1\cos\big(\frac{\a}{2}(\omega-\pi)-\omega\big)
+\lambda_0\cos\big(\frac{\a}{2}(\omega-\pi)\big)+\lambda_{-1}\cos\big(\frac{\a}{2}(\omega-\pi)+\omega\big),
\end{equation}
where $\omega\in[0,\pi]$ and $\lambda_1, \lambda_0, \lambda_{-1}$ are defined in \eqref{eq:c1}.
Then $h(\a,\omega)$ does not decrease with respect to $\omega$.
\end{Lem}

The proof of above lemma is elementary but quite technical, and hence
deferred to the appendix.
%
\begin{Lem}\label{lem:une}
For $1<\a\leqslant 2$, we have
\begin{equation}\label{eq:une}
C_{\a}|u|_{H^{\a/2}_h}^2\leqslant
({\Delta}^{\a}_{h}u,u)_h\leqslant
|u|_{H^{\a/2}_h}^2,
\end{equation}
where $C_{\a}=\frac{2^{\a}(1-\a^2)}{3\pi^{\a}\cos{\frac{\a\pi}{2}}}>0$.
\end{Lem}
\textbf{Proof.} From the Parseval's identity \eqref{eq:PI}, it follows that
\begin{equation}\label{eq:relat1}
({\Delta}^{\a}_{h}u,u)_h=\int^{+\pi/h}_{-\pi/h}h^{-\a}f(\a,k)\widehat{u}(k)\overline{\widehat{u}(k)}\text{d}k,
\end{equation}
where
\begin{equation}
f(\a,k)=\frac{1}{2\cos{\frac{\a\pi}{2}}}\Big(\sum^{\infty}_{j=0}w^{(\a)}_je^{i(j-1)hk}+\sum^{\infty}_{j=0}w^{(\a)}_je^{-i(j-1)hk}\Big).
\end{equation}
In view of \eqref{eq:cpq} and \eqref{eq:GL}, we get
\begin{eqnarray*}
f(\a,k)&=&\frac{1}{2\cos{\frac{\a\pi}{2}}}\Big[\lambda_1e^{-ihk}\sum^{\infty}_{j=0}g^{(\a)}_je^{ijhk}+\lambda_0\sum^{\infty}_{j=0}g^{(\a)}_je^{ijhk}
+\lambda_{-1}e^{ihk}\sum^{\infty}_{j=0}g^{(\a)}_je^{ijhk}\\
&{}&+\lambda_1e^{ihk}\sum^{\infty}_{j=0}g^{(\a)}_je^{-ijhk}+\lambda_0\sum^{\infty}_{j=0}g^{(\a)}_je^{-ijhk}
+\lambda_{-1}e^{-ihk}\sum^{\infty}_{j=0}g^{(\a)}_je^{-ijhk}\Big]\\
&=&\frac{1}{2\cos{\frac{\a\pi}{2}}}\Big[\lambda_1\Big(e^{-ihk}(1-e^{ihk})^{\a}+e^{ik}(1-e^{-ihk})^{\a}\Big)
+\lambda_0\Big((1-e^{ihk})^{\a}+(1-e^{-ihk})^{\a}\Big)\\
&{}&+\lambda_{-1}\Big(e^{ihk}(1-e^{ihk})^{\a}+e^{-ik}(1-e^{-ihk})^{\a}\Big)\Big].
\end{eqnarray*}
Clearly, $f(\a,k)$ is a real-valued even function, and it is
therefore sufficient to consider its principle value for
$k\in[0,\pi/h]$. Invoking the relation
$e^{i\theta}-e^{i\phi}=2i\sin\big(\frac{\theta-\phi}{2}\big)e^{\frac{i(\theta+\phi)}{2}}$,
 we have
\begin{equation*}
\begin{split}
f(\a,k)
&=\frac{\Big(2\sin\frac{hk}{2}\Big)^{\a}}{\cos{\frac{\a\pi}{2}}}\Big(\lambda_1\cos\big(\frac{\a}{2}(hk-\pi)-hk\big)
+\lambda_0\cos\big(\frac{\a}{2}(hk-\pi)\big)+\lambda_{-1}\cos\big(\frac{\a}{2}(hk-\pi)+hk\big)\Big)\\
&=\frac{\Big(2\sin\frac{hk}{2}\Big)^{\a}}{\cos{\frac{\a\pi}{2}}}h(\a,hk),
\end{split}
\end{equation*}
where $h(\a,hk)$ is defined as in \eqref{eq:hak} with $\omega=hk\in[0,\pi]$ .
From Lemma \ref{Lem:L1}, it follows that
\begin{equation}\label{eq:cosbound1}
\cos\big(\frac{\a\pi}{2}\big)=h(\a,0)\leqslant h(\a,hk)\leqslant h(\a,\pi)=\frac{1-\a^2}{3}.
\end{equation}
Hence, combining above inequality with the fact that for $k\in[0,\pi/h]$, $hk/\pi\leqslant\sin(hk/2)\leqslant hk/2$, we obtain
\begin{equation}
C_{\a}|hk|^{\a}\leqslant f(\a,k)\leqslant |hk|^{\a}.
\end{equation}
This together with \eqref{eq:relat1} implies \eqref{eq:une} and thus completes the proof.
$\Box$
\begin{Rem}
This idea can be used to analyze some other popular second order
schemes, including the WSGD methods proposed in \cite{TianW:2015:MC}
and the fractional centered difference method
\cite{Ortigueira:2006:IJMMS,CelikC:2012:JCP}, and similar results
can be derived. Unfortunately, it seems difficult to extend directly
this study to fourth order schemes, such as the weighted and shifted
Lubich difference method \cite{ChenM:2014:SIAMNUA,ChenM:2014:CICP},
because the involved function $f(\a,k)$ is much more complicated.
\end{Rem}
\section{Finite difference scheme}
%
In practical computation, the
whole space problem is usually truncated onto a finite interval
$\Omega=(a,b)$ subject to the homogeneous boundary condition ($a$ and $b$ are usually chosen sufficient large such
that the truncation error is negligible). Thus the FGLE \eqref{eq:FSE1}-\eqref{eq:IC1} is truncated on the interval $\Omega=(a,b)$ as
\begin{eqnarray}
&&u_t+(\upsilon+i\eta)(-\Delta)^{\frac{\a}{2}}u+(\kappa+i\zeta){|u|}^2u-\gamma u=0,\quad x\in\Omega,\ t\in(0,T],\label{eq:FSE}\\
&&u(x,0)=u_0(x), \quad x\in\mathbb{R},\label{eq:IC}\\
&&u(x,t)=0, \quad x\in \mathbb{R}\backslash\Omega,\ t\in[0,T]. \label{eq:BC}
\end{eqnarray}
The boundary condition \eqref{eq:BC} is referred to as the nonlocal volume constraint (or the extended Dirichlet boundary) and the corresponding problem \eqref{eq:FSE}-\eqref{eq:BC}
as the volume constraint problem (see \cite{DuQ:2012:SIAMR,DefterliO:2015:FCAA} for more details). It is noted that, under this boundary, the fractional derivative has reduced to $-(-\Delta)^{\frac{\a}{2}}u(x,t)=-\frac{1}{2\cos{\frac{\alpha\pi}{2}}}\big[{}_{a}D^{\alpha}_{x}u(x,t)+{}_{x}D^{\alpha}_{b}u(x,t)\big]$.

Let $h=\frac{b-a}{M}$ with a positive integer $M$ and define $
x_j=a+jh,\ 0\leqslant j \leqslant M$. Owing to the above boundary
constraint \eqref{eq:BC}, if $u\in\mathscr{L}^{2+\a}(\mathbb{R})$ (see Remark 2.5 in \cite{HaoZ:2015:JCP}), the WSGD operator \eqref{eq:WSGDR1} can be
simplified as
\begin{equation}\label{eq:WSGDR2}
\begin{split}
&{}_{L}\mathcal{D}^{\a}_{h}u(x_j):=\frac{1}{h^{\a}}\sum^{j+1}_{l=0}w^{(\a)}_lu(x_{j-l+1})={}_{a}D^{\a}_{x}u(x_j)+O(h^2),\\
&{}_{R}\mathcal{D}^{\a}_{h}u(x_j):=\frac{1}{h^{\a}}\sum^{M-j+1}_{l=0}w^{(\a)}_lu(x_{j+l-1})={}_{x}D^{\a}_{b}u(x_j)+O(h^2),
\end{split}
\end{equation}
and WSGD approximation \eqref{eq:wsgdr1} for the fractional Laplacian as
\begin{equation}\label{eq:wsgdr2}
\begin{split}
\Delta^{\a}_{h}u(x_j)&=\frac{1}{h^{\a}}\frac{1}{2\cos{\frac{\a\pi}{2}}}\Big(\sum^{j+1}_{l=0}w^{(\a)}_lu(x_{j-l+1})+\sum^{M-j+1}_{l=0}w^{(\a)}_lu(x_{j+l-1})\Big)\\
&=(-\Delta)^{\frac{\a}{2}}u(x_j)+O(h^2).
\end{split}
\end{equation}

Chosen the time step $\tau:=\frac{T}{N}$ with a positive integer
$N$,
define a partition of $[0,T]\times[a,b]$ by
${\Omega}_{\tau}\times {\Omega}_{h}$ with the grid
${\Omega}_{\tau}=\{t_n\ |\ t_n=n\tau,\ 0\leqslant n\leqslant N\}$
and ${\Omega}_{h}=\{x_j\ |\ x_j=a+jh,\ 0\leqslant j \leqslant M\}$.

Given a grid function $v=\{v^n_j|(x_j,t_n)\in{\Omega}_{\tau}\times
{\Omega}_{h}\}$, denote
\begin{equation}
\delta_{t}v^{n+\frac{1}{2}}_j=\frac{v^{n+1}_{j}-v^{n}_{j}}{\tau},\quad
v^{n+\frac{1}{2}}_j=\frac{v^{n+1}_{j}+v^{n}_{j}}{2}.
\end{equation}
Denote the index set $\mathcal{T}_M=\{j\ |\ j=1,2,\ldots,M-1\}$ and the grid function space $\mathcal{V}_h=\{v\ | \
v=(v_1,v_2,\ldots,v_{M-1})\}$.

Under the boundary constraint \eqref{eq:BC}, the inner product $(\cdot,\cdot)_h$ and norms $\|\cdot\|_h$, $\|\cdot\|_{l^p_h}$, $\|\cdot\|_{l^{\infty}_h}$ previously defined in the unbounded
interval carry over to the finite interval by regarding that $u_j=0$ for $j\leqslant 0$ and $j\geqslant M$. Hence in these notations, we just restrict the index $j$ from $1$ to $M-1$ and continue to
use these notations without confusion for convenience. Based on these considerations, the inequalities introduced in above section still hold in the finite interval.

With these premises, we now propose a difference scheme for the FGLE
\eqref{eq:FSE}. Let $u^n_j$ be the numerical approximation of
$u(x_j,t_n)$. Applying the implicit midpoint method in time and the
WSGD approximation for the fractional Laplacian, the difference
scheme reads
\begin{eqnarray}
&&\delta_{t}u^{n+\frac{1}{2}}_j+(\upsilon+i\eta)\Delta^{\a}_{h}u^{n+\frac{1}{2}}_j+(\kappa+i\zeta)|u^{n+\frac{1}{2}}_j|^2u^{n+\frac{1}{2}}_j-\gamma u^{n+\frac{1}{2}}_j=0,\label{eq:CCNS}\\
&&\ \ \ \ \ \ \ \ \ \ \ \ \ \ \ \ \ \ \ \ \ \ \ \ \ \ \ \ \ \ \ \  \ \ \ \ \ \ \ \ \ \ \ \ \ \ \ \ \ \ \ \ \ \ \ \ \ \ \ \quad j\in\mathcal{T}_M, \ 0\leqslant n \leqslant N-1,\nonumber\\
&&u^0_j=u_0(x_j), \quad j\in\mathbb{Z},\label{eq:DIC}\\
&&u^n_j=0, \quad j\in\mathbb{Z}\backslash\mathcal{T}_M, \quad 0\leqslant n \leqslant N.\label{eq:DBC}
\end{eqnarray}
%
\section{Theoretical analysis}
%
In this section, we study theoretical properties of the scheme
\eqref{eq:CCNS}-\eqref{eq:DBC}, including the a priori estimate,
solvability and convergence.
\subsection{A priori bound}
For showing the a priori bound of the solution to the scheme \eqref{eq:CCNS}-\eqref{eq:DBC}, we first introduce some notations and lemmas.

Denote matrix
\[
 \mathbf{W}=
 \left( \begin{array}{ccccc}
 w^{(\a)}_1     & w^{(\a)}_0     & {}         & {}         & {} \\
 w^{(\a)}_2     & w^{(\a)}_1     & w^{(\a)}_0 & {}         & {} \\
 \vdots         & w^{(\a)}_2     & w^{(\a)}_1 & \ddots     & {} \\
 w^{(\a)}_{M-2} & \vdots         & \ddots     & \ddots     & w^{(\a)}_0\\
 w^{(\a)}_{M-1} & w^{(\a)}_{M-2} & \cdots     & w^{(\a)}_2 & w^{(\a)}_1
 \end{array} \right)\in\mathbb{R}^{(M-1)\times(M-1)},
\]
and matrix
$\mathbf{C}=\frac{1}{2\cos{\frac{\a\pi}{2}}}(\mathbf{W}+{\mathbf{W}}^{T})\in\mathbb{R}^{(M-1)\times(M-1)}$.
Then for $u=(u_1,u_2,\ldots,u_{M-1})^{T}$, we can rewrite
$\Delta^{\a}_{h}u=\frac{1}{h^{\a}}\mathbf{C}u$. In addition, we have
the following lemmas.
\begin{Lem}
Matrix $\mathbf{C}$ is a real-valued symmetry positive definite matrix.
\end{Lem}
\textbf{Proof.} It is obviously seen that $\mathbf{C}$ is a real-valued symmetry matrix. The positive definiteness, for $1<\a<2$, can be obtained by invoking the property of coefficients in \eqref{eq:cpq1} and the Ger\v{s}gorin disc theorem \cite{HornRA:1986:CUP}. For $\a=2$, matrix $\mathbf{C}$ reduces to the classical Laplacian matrix associated with the second order centered difference and thus, the positive definiteness is obtained.
$\Box$
%
\begin{Lem}\label{lem:BG}
For any two grid functions $u,v\in{\mathcal{V}}_h$, there exists a
linear operator $\Lambda^\a$ such that
 \begin{equation}\label{eq:RELATION}
   \Big(\Delta^{\a}_{h}u,v\Big)_h=\Big(\Lambda^\a u,\Lambda^\a v\Big)_h.
 \end{equation}
\end{Lem}
\textbf{Proof.} The proof is similar to that in \cite{WangP:2015:JCP} (see Lemma 3.1 in \cite{WangP:2015:JCP})
where the fractional centered difference is adopted.
The linear operator $\Lambda^\a$ is defined by $\Lambda^\a u=h^{-\frac{\a}{2}}\mathbf{L} u$ where matrix $\mathbf{L}$ satisfying $\mathbf{C}=\mathbf{L}^T\mathbf{L}$ is the Cholesky factor.
$\Box$

Based on the previous lemma, we can establish the following boundedness
estimate.
\begin{Theo}\label{th:bound}
The difference solution of scheme \eqref{eq:CCNS}-\eqref{eq:DBC} is
bounded in the following sense
\begin{equation}\label{eq:bound}
\|u^n\|_h\leqslant C_M,\quad 0\leqslant n \leqslant N.
\end{equation}
\end{Theo}
\textbf{Proof.} Computing the discrete inner product of
\eqref{eq:CCNS} with $u^{n+\frac{1}{2}}$, then taking the real
part of the resulting equation, we obtain
\begin{equation}\label{eq:thbound1}
\frac{\|u^{n+1}\|^2_h-\|u^{n}\|^2_h}{2\t}+\upsilon\|\Lambda^\a u^{n+\frac{1}{2}}\|^2_h+\kappa\|u^{n+\frac{1}{2}}\|^4_{l^4_h}=\gamma\|u^{n+\frac{1}{2}}\|^2_h, \quad 0\leqslant n \leqslant N-1,
\end{equation}
where we have used the relation \eqref{eq:RELATION}.

If $\gamma\leqslant 0$, from the above inequality we get
\[
\|u^{n+1}\|^2_h\leqslant\|u^{n}\|^2_h\leqslant\cdots\leqslant\|u^{0}\|^2_h,\quad 0\leqslant n \leqslant N-1.
\]

If $\gamma> 0$, it follows from \eqref{eq:thbound1} that
\[
\|u^{n+1}\|^2_h-\|u^{n}\|^2_h\leqslant 2\t\gamma\|u^{n+\frac{1}{2}}\|^2_h\leqslant \t\gamma\big(\|u^{n+1}\|^2_h+\|u^{n}\|^2_h\big),\quad 0\leqslant n \leqslant N-1.
\]
Let $\t\leqslant\frac{1}{2\gamma}$. We have
\[
\|u^{n+1}\|^2_h\leqslant\big(1+4\gamma\t\big)\|u^{n}\|^2_h,\quad 0\leqslant n \leqslant N-1,
\]
which immediately implies
\[
\|u^{n}\|^2_h\leqslant \exp(4\gamma T)\|u^{0}\|^2_h,\quad 0\leqslant n \leqslant N.
\]
Thus the proof is complete.
$\Box$
\begin{Rem}
Theorem \ref{th:bound} implies that the difference solution of
scheme \eqref{eq:CCNS}-\eqref{eq:DBC} is bounded for a long time if
$\gamma\leqslant 0$, while is bounded for a given $T$ if $\gamma>
0$. 
\end{Rem}
\subsection{Solvability}
The existence of the solution is shown by virtue of the Brouwder fixed point theorem.

\begin{Lem}(Brouwder fixed point theorem \cite{AkrivisG:1991:NM})\label{Lem:EXIST}
Let $(\mathcal{H},\langle\cdot,\cdot\rangle)$ be a finite
dimensional inner product space, $\|\cdot\|$ be the associated norm,
and $f:\mathcal{H}\rightarrow \mathcal{H}$ be continuous. Assume,
moreover, that
\begin{equation}
 \exists \rho>0,\ \forall z\in \mathcal{H},\ \|z\|=\rho,\ Re\langle g(z),\ z\rangle \geqslant 0.
\end{equation}
Then, there exists a $z^{*}\in \mathcal{H}$ such that $g(z^{*})=0$
and $\|z^{*}\|\leqslant\rho.$
\end{Lem}
\begin{Theo}
The solution of difference scheme \eqref{eq:CCNS}-\eqref{eq:DBC} exists.
\end{Theo}
\textbf{Proof.} The proof proceeds in an inductive way. Obviously,
$u^0$ has been determined uniquely from \eqref{eq:DIC} and
\eqref{eq:DBC}. For given $u^n$ ($0\leqslant n\leqslant N-1$), it
remains to prove that there exists $u^{n+1}$ satisfying the scheme.
To this end, for fixed $n$, rewrite \eqref{eq:CCNS} in the form
\begin{equation*}
u^{n+\frac{1}{2}}_j=u^n_j-\frac{\t}{2}\Big[(\upsilon+i\eta)\Delta^{\a}_{h}u^{n+\frac{1}{2}}_j+(\kappa+i\zeta)|u^{n+\frac{1}{2}}_j|^2u^{n+\frac{1}{2}}_j-\gamma u^{n+\frac{1}{2}}_j\Big], \quad j\in\mathcal{T}_M.
\end{equation*}
Consider a mapping $\mathcal{F}: \mathcal{V}_h\rightarrow\mathcal{V}_h$ defined by
\begin{equation}\label{eq:ext1}
\big(\mathcal{F}(v)\big)_j=v_j-u^n_j+\frac{\t}{2}\Big[(\upsilon+i\eta)\Delta^{\a}_{h}v_j+(\kappa+i\zeta)|v_j|^2v_j-\gamma v_j\Big],\quad j\in\mathcal{T}_M,
\end{equation}
which is obviously continuous. Computing the discrete inner product of \eqref{eq:ext1} with $v$ gives
\begin{equation}
\begin{split}
\big(\mathcal{F}(v),v\big)_h
&=\|v\|^2_h-(u^n,v)_h+\frac{\t}{2}\Big[(\upsilon+i\eta)\|\Lambda^\a v\|^2_h+(\kappa+i\zeta)\|v\|^4_{l^4_h}-\gamma \|v\|^2_h\Big],
\end{split}
\end{equation}
where \eqref{eq:RELATION} was used. Then, taking the real part, we obtain
\begin{equation}
\begin{split}
Re\big(\mathcal{F}(v),v\big)_h&=(1-\frac{\gamma\t}{2})\|v\|^2_h-Re(u^n,v)+\frac{\upsilon\t}{2}\|\Lambda^\a v\|^2_h+\frac{\kappa\t}{2}\|v\|^4_{l^4_h}\\
&\geqslant (1-\frac{\gamma\t}{2})\|v\|^2_h-\|u^n\|_h\|v\|_h\\
&=\|v\|_h\big((1-\frac{\gamma\t}{2})\|v\|_h-\|u^n\|_h\big).
\end{split}
\end{equation}
Hence, taking $\gamma\t\leqslant 1$ and $\|v\|_h=2\|u^n\|_h$, there exists $Re\big(\mathcal{F}(v),v\big)_h\geqslant 0$. Then from Lemma \ref{Lem:EXIST}, we obtain the existence of $u^{n+\frac{1}{2}}$ and thus, the existence of $u^{n+1}$ by noting that $u^{n+1}=2u^{n+\frac{1}{2}}-u^n$.
$\Box$

For a technical reason, the uniqueness of the solution will be shown in the subsection 4.4.
\subsection{Convergence}
Before establishing the convergence, we first analyze the local
truncation error of scheme \eqref{eq:CCNS}-\eqref{eq:DBC}. For notational convenience we denote
grid functions $U^n_j:=u(x_j,t_n)$. Define the truncation error as
\begin{equation}\label{eq:err}
\begin{split}
  R^{n+\frac{1}{2}}_j:=\delta_{t}U^{n+\frac{1}{2}}_j+(\upsilon+i\eta)\Delta^{\a}_{h}U^{n+\frac{1}{2}}_j+(\kappa+i\zeta)|U^{n+\frac{1}{2}}_j|^2U^{n+\frac{1}{2}}_j-\gamma U^{n+\frac{1}{2}}_j,&\\
  j\in\mathcal{T}_M, \quad 0\leqslant n \leqslant N-1.&
\end{split}
\end{equation}
Then from \eqref{eq:wsgdr2} and Taylor's expansion, we can obtain the following local truncation error estimate.
\begin{Lem}\label{lem:locerr}
Suppose that the problem \eqref{eq:FSE}-\eqref{eq:BC} has a smooth solution. Then we have
\begin{equation}\label{eq:Err1}
|R^{n+\frac{1}{2}}_j|\leqslant C_R({\tau}^2+h^2), \quad
j\in\mathcal{T}_M,\ 0\leqslant n \leqslant N-1.
\end{equation}
\end{Lem}
\textbf{Proof.}
Applying the Taylor's expansion of the solution at $(x_j,t_{n+\frac{1}{2}})$ yields
\begin{equation}\label{eq:lte1}
\begin{split}
&\frac{u(x_j,t_{n+1})-u(x_j,t_n)}{\t}\\
&=\partial_t u\big(x_j,t_{n+\frac{1}{2}}\big)+\frac{\t^2}{16}\int^1_0\Big[\frac{\partial^3u}{\partial t^3}\big(x_j,t_{n+\frac{1}{2}}+\frac{s}{2}\t\big)+\frac{\partial^3u}{\partial t^3}\big(x_j,t_{n+\frac{1}{2}}-\frac{s}{2}\t\big)\Big](1-s)^2\text{d}s,
\end{split}
\end{equation}
and
\begin{equation}\label{eq:lte2}
\begin{split}
&\frac{u(x_j,t_{n+1})+u(x_j,t_n)}{2}\\
&=u\big(x_j,t_{n+\frac{1}{2}}\big)+\frac{\t^2}{8}\int^1_0\Big[\frac{\partial^2u}{\partial t^2}\big(x_j,t_{n+\frac{1}{2}}+\frac{s}{2}\t\big)+\frac{\partial^2u}{\partial t^2}\big(x_j,t_{n+\frac{1}{2}}-\frac{s}{2}\t\big)\Big](1-s)\text{d}s.
\end{split}
\end{equation}
Furthermore, noticing the error estimate \eqref{eq:wsgdr2}, we have
\begin{equation}\label{eq:lte3}
\begin{split}
&\Delta^{\a}_{h}\Big(\frac{u(x_j,t_{n+1})+u(x_j,t_n)}{2}\Big)-(-\Delta)^{\frac{\a}{2}} u\big(x_j,t_{n+\frac{1}{2}}\big)\\
&=\frac{\t^2}{8}\int^1_0(-\Delta)^{\frac{\a}{2}}\Big[\frac{\partial^2u}{\partial t^2}\big(x_j,t_{n+\frac{1}{2}}+\frac{s}{2}\t\big)+\frac{\partial^2u}{\partial t^2}\big(x_j,t_{n+\frac{1}{2}}-\frac{s}{2}\t\big)\Big](1-s)\text{d}s+O(h^2),
\end{split}
\end{equation}
Substituting \eqref{eq:lte1}-\eqref{eq:lte3} into \eqref{eq:err} gives
\eqref{eq:Err1} immediately and thus,
completes the proof.
$\Box$

From above lemma, it follows  that
\begin{equation}\label{eq:Err11}
\|R^{n+\frac{1}{2}}\|^2_h\leqslant
(b-a)\big(C_R({\tau}^2+h^2)\big)^2, \quad 0\leqslant n \leqslant
N-1.
\end{equation}
Define the error function $e^n\in\mathcal {V}_h$ as
\[
 e^n_j=U^n_j-u^n_j, \quad j\in\mathcal{T}_M,\quad 0\leqslant n
\leqslant N.
\]
Then we get the following convergence result.
\begin{Theo}\label{th:Converg}
Suppose that the problem \eqref{eq:FSE}-\eqref{eq:BC} has a smooth solution. Then there exists ${\tau}_0>0$
sufficiently small such that, when $0<\tau\leqslant{\tau}_0$, we have
\begin{equation}\label{eq:aabound}
\|e^n\|_h \leqslant C(\tau^2+h^2), \quad 0\leqslant n \leqslant N,
\end{equation}
where $C$ denotes a positive constant independent of $\t$ and $h$.
\end{Theo}
\textbf{Proof.}
Subtracting \eqref{eq:CCNS} from \eqref{eq:err} gives the following
error equation
\begin{equation}\label{eq:ee1}
\begin{split}
&\delta_{t}e^{n+\frac{1}{2}}_j+(\upsilon+i\eta)\Delta^{\a}_{h}e^{n+\frac{1}{2}}_j+(\kappa+i\zeta)|u^{n+\frac{1}{2}}_j|^2e^{n+\frac{1}{2}}_j
+(\kappa+i\zeta)G^{n+\frac{1}{2}}_jU^{n+\frac{1}{2}}_j-\gamma e^{n+\frac{1}{2}}_j=R^{n+\frac{1}{2}}_j,\\
&\ \ \ \ \ \ \ \ \ \ \ \ \ \ \ \ \ \ \ \ \ \ \ \ \ \ \ \ \ \ \ \ \ \ \ \ \ \ \ \ \ \ \ \ \ \ \ \ \ \ \ \ \ \ \ \ \ \ \ \ \ \  \ \ \ \ \ \ \ \ \ \ \ \ \ \ \ \ \ \ \ \ \ \ \ \ \ \ \ \ \ \ \ \ \ \  \ \quad j\in\mathcal{T}_M,\quad 0\leqslant n \leqslant N-1,
\end{split}
\end{equation}
where
\begin{equation*}
G^{n+\frac{1}{2}}_j=|U^{n+\frac{1}{2}}_j|^2-|u^{n+\frac{1}{2}}_j|^2.
\end{equation*}
Computing the discrete inner product of
\eqref{eq:ee1} with $e^{n+\frac{1}{2}}$, then taking the real
part of the resulting equation, we obtain
\begin{equation}\label{eq:ee2}
\begin{split}
&\frac{\|e^{n+1}\|^2_h-\|e^{n}\|^2_h}{2\t}+\upsilon\|\Lambda^\a e^{n+\frac{1}{2}}\|^2_h+\kappa h\sum^{M-1}_{j=1}|u^{n+\frac{1}{2}}_j|^2|e^{n+\frac{1}{2}}_j|^2-\gamma\|e^{n+\frac{1}{2}}\|^2_h\\
&+Re\Big((\kappa+i\zeta)h\sum^{M-1}_{j=1}G^{n+\frac{1}{2}}_jU^{n+\frac{1}{2}}_j\overline{e^{n+\frac{1}{2}}_j}\Big)
=Re(R^{n+\frac{1}{2}},e^{n+\frac{1}{2}})_h, \quad 0\leqslant n \leqslant N-1.
\end{split}
\end{equation}
We first estimate the last term on the left-hand side of
\eqref{eq:ee2}. In view of the smoothness assumption of the exact
solution, denote $C_u=\sup_{0\leqslant t \leqslant
T, a\leqslant x \leqslant b}|u(x,t)|$. Noting that
\[
|u^n_j|\leqslant|U^n_j|+|e^n_j|\leqslant C_u+|e^n_j|,\quad 0\leqslant n \leqslant N,
\]
we have for $0\leqslant n \leqslant N-1$,
\[
\begin{split}
\Big|G^{n+\frac{1}{2}}_jU^{n+\frac{1}{2}}_j\Big|&=\Big|\big(|U^{n+\frac{1}{2}}_j|^2-|u^{n+\frac{1}{2}}_j|^2\big)U^{n+\frac{1}{2}}_j\Big|\\
&\leqslant C_u|U^{n+\frac{1}{2}}_j-u^{n+\frac{1}{2}}_j|\big(|U^{n+\frac{1}{2}}_j|+|u^{n+\frac{1}{2}}_j|\big)\\
&\leqslant C_u|e^{n+\frac{1}{2}}_j|\big(2C_u+|e^{n+\frac{1}{2}}_j|\big)\\
&\leqslant 2C_u^2|e^{n+\frac{1}{2}}_j|+C_u|e^{n+\frac{1}{2}}_j|^2.
\end{split}
\]
Then we obtain
\begin{equation}\label{eq:G1}
\begin{split}
&\Big|Re\Big((\kappa+i\zeta)h\sum^{M-1}_{j=1}G^{n+\frac{1}{2}}_jU^{n+\frac{1}{2}}_j\overline{e^{n+\frac{1}{2}}_j}\Big)\Big|\\
&\leqslant \sqrt{\kappa^2+\zeta^2}h\sum^{M-1}_{j=1}(2C_u^2|e^{n+\frac{1}{2}}_j|+C_u|e^{n+\frac{1}{2}}_j|^2)|e^{n+\frac{1}{2}}_j|\\
&= \sqrt{\kappa^2+\zeta^2}(2C_u^2\|e^{n+\frac{1}{2}}\|^2_h+C_u\|e^{n+\frac{1}{2}}\|^3_{l^3_h}), \quad 0\leqslant n \leqslant N-1.
\end{split}
\end{equation}
For the term on the right-hand side of \eqref{eq:ee2}, using the
Cauchy-Swcharz inequality gives
\begin{equation}\label{eq:R1}
Re(R^{n+\frac{1}{2}},e^{n+\frac{1}{2}})_h\leqslant \frac{1}{2}(\|R^{n+\frac{1}{2}}\|^2_h+\|e^{n+\frac{1}{2}}\|^2_h), \quad 0\leqslant n \leqslant N-1.
\end{equation}
By substituting \eqref{eq:G1} and \eqref{eq:R1} into \eqref{eq:ee2}, we get for $0\leqslant n \leqslant N-1$,
\begin{equation}\label{eq:ee3}
\begin{split}
&\frac{\|e^{n+1}\|^2_h-\|e^{n}\|^2_h}{2\t}+\upsilon\|\Lambda^\a e^{n+\frac{1}{2}}\|^2_h\\
&\leqslant  \sqrt{\kappa^2+\zeta^2}(2C_u^2\|e^{n+\frac{1}{2}}\|^2_h+C_u\|e^{n+\frac{1}{2}}\|^3_{l^3_h})
+|\gamma|\|e^{n+\frac{1}{2}}\|^2_h+\frac{1}{2}(\|R^{n+\frac{1}{2}}\|^2_h+\|e^{n+\frac{1}{2}}\|^2_h)\\
&=\big(2\sqrt{\kappa^2+\zeta^2}C_u^2+|\gamma|+\frac{1}{2}\big)\|e^{n+\frac{1}{2}}\|^2_h+\sqrt{\kappa^2+\zeta^2}C_u\|e^{n+\frac{1}{2}}\|^3_{l^3_h}
+\frac{1}{2}\|R^{n+\frac{1}{2}}\|^2_h.
\end{split}
\end{equation}
Furthermore, Theorem \ref{th:bound} implies that
\begin{equation}
\|e^{n+\frac{1}{2}}\|_h\leqslant \|U^{n+\frac{1}{2}}\|_h+\|u^{n+\frac{1}{2}}\|_h\leqslant \sqrt{b-a}C_u+C_M, \quad 0\leqslant n \leqslant N-1.
\end{equation}
Hence, in view of \eqref{eq:dgni} with $p=3$ and
$\sigma_0=\frac{\a}{6}$, we have
\begin{equation}
\begin{split}
\|e^{n+\frac{1}{2}}\|^3_{l^3_h}&\leqslant C_{\sigma_0}\|e^{n+\frac{1}{2}}\|^2_h\|e^{n+\frac{1}{2}}\|_{H^{\a/2}_h}\\
&\leqslant C_{\sigma_0}(\varepsilon|e^{n+\frac{1}{2}}|_{H^{\a/2}_h}^2+\varepsilon\|e^{n+\frac{1}{2}}\|^2_h+\frac{1}{4\varepsilon}\|e^{n+\frac{1}{2}}\|^4_h)\\
&\leqslant C_{\sigma_0}\Big(\frac{\varepsilon}{C_{\a}}\|\Lambda^\a e^{n+\frac{1}{2}}\|^2_h+\big(\varepsilon+\frac{(\sqrt{b-a}C_u+C_M)^2}{4\varepsilon}\big)\|e^{n+\frac{1}{2}}\|^2_h\Big),
\end{split}
\end{equation}
where we have used \eqref{eq:une} and \eqref{eq:RELATION} for the last inequality. Taking $\varepsilon=\frac{C_{\a}\upsilon}{\sqrt{\kappa^2+\zeta^2}C_uC_{\sigma_0}}$ and plugging the above inequality into \eqref{eq:ee3}, we obtain
\begin{equation*}\label{eq:ee4}
\frac{\|e^{n+1}\|^2_h-\|e^{n}\|^2_h}{2\t}+\upsilon\|\Lambda^\a e^{n+\frac{1}{2}}\|^2_h
\leqslant \upsilon\|\Lambda^\a e^{n+\frac{1}{2}}\|^2_h+C_1\|e^{n+\frac{1}{2}}\|^2_h+\frac{1}{2}\|R^{n+\frac{1}{2}}\|^2_h, \quad 0\leqslant n \leqslant N-1,
\end{equation*}
namely,
\begin{equation}
\|e^{n+1}\|^2_h-\|e^{n}\|^2_h\leqslant \t C_1(\|e^{n+1}\|^2_h+\|e^{n}\|^2_h)+\t \|R^{n+\frac{1}{2}}\|^2_h, \quad 0\leqslant n \leqslant N-1,
\end{equation}
where $C_1=2\sqrt{\kappa^2+\zeta^2}C_u^2+|\gamma|+\frac{1}{2}+C_{\a}\upsilon+\frac{(\kappa^2+\zeta^2)C_u^2C_{\sigma_0}^2(\sqrt{b-a}C_u+C_M)^2}{4C_{\a}\upsilon}$.

If $\t\leqslant \frac{1}{2C_1}$, we have
\begin{equation}
\|e^{n+1}\|^2_h\leqslant (1+4C_1\t)\|e^{n}\|^2_h+2\t \|R^{n+\frac{1}{2}}\|^2_h, \quad 0\leqslant n \leqslant N-1.
\end{equation}
This together with \eqref{eq:Err11} and the discrete Gronwall inequality gives
\begin{equation}
\|e^{n}\|^2_h\leqslant \exp(4C_1T)\frac{(b-a)C_R^2}{2C_1}(\t^2+h^2)^2, \quad 0\leqslant n \leqslant N,
\end{equation}
which implies \eqref{eq:aabound} with $C=\exp(2C_1T)C_R\sqrt{\frac{b-a}{2C_1}}$. Thus, the proof is complete.
$\Box$
\subsection{Uniqueness}
Now we are in a position to show the uniqueness of the solution to the scheme \eqref{eq:CCNS}-\eqref{eq:DBC}.

From Theorem \ref{th:Converg}, using the inverse inequality $\|\cdot\|_{l^{\infty}_h}^2\leqslant h^{-1}\|\cdot\|^2_h$, we first show the uniform boundedness of the difference solution. In fact, assume $\t\leqslant C_2h$, then for $0<h\leqslant (\frac{1}{C(1+C_2^2)})^{\frac{2}{3}}$, we obtain
\begin{equation}\label{eq:unib1}
\|u^n\|_{l^{\infty}_h}\leqslant \|U^n\|_{l^{\infty}_h}+\|e^n\|_{l^{\infty}_h}\leqslant C_u+Ch^{-\frac{1}{2}}(\t^2+h^2)\leqslant 1+C_u, \quad 0\leqslant n \leqslant N.
\end{equation}
Using this inequality, we have the following results.
\begin{Theo}
Assume $\t\leqslant Ch$, then there exist $\t_0>0$ and $h_0>0$ sufficiently small such that, when $0<\t\leqslant \t_0$ and $0<h\leqslant h_0$,  the solution of difference scheme \eqref{eq:CCNS}-\eqref{eq:DBC} is unique.
\end{Theo}
\textbf{Proof.} Suppose there exist two solutions $u^{(1)}, u^{(2)}\in\mathcal{V}_h$ to the scheme \eqref{eq:CCNS}-\eqref{eq:DBC}. Then from \eqref{eq:unib1}, we have
\begin{equation}\label{eq:uni1}
\|u^n\|_{l^{\infty}_h}\leqslant 1+C_u, \quad \|u^{(1)}\|_{l^{\infty}_h}\leqslant 1+C_u,\quad \|u^{(2)}\|_{l^{\infty}_h}\leqslant 1+C_u, \quad 0\leqslant n \leqslant N-1.
\end{equation}
Setting $w=u^{(1)}-u^{(2)}$, we obtain
\begin{equation}\label{eq:uni2}
\frac{w_j}{\t}+\frac{1}{2}(\upsilon+i\eta)\Delta^{\a}_{h}w_j+(\kappa+i\zeta)g_j-\frac{1}{2}\gamma w_j=0,\quad j\in\mathcal{T}_M,
\end{equation}
where
\[
g_j=\Big|\frac{u^{(1)}_j+u^n_j}{2}\Big|^2\frac{u^{(1)}_j+u^n_j}{2}-\Big|\frac{u^{(2)}_j+u^n_j}{2}\Big|^2\frac{u^{(2)}_j+u^n_j}{2}.
\]
Computing the discrete inner product of \eqref{eq:uni2} with $w$ and taking the real part give
\begin{equation}\label{eq:uni3}
\|w\|^2_h+\frac{\upsilon\t}{2}\|\Lambda^\a w\|^2_h+\t Re(\kappa+i\zeta)(g,w)_h=\frac{\gamma\t}{2}\|w\|^2_h.
\end{equation}
For the third term on the left-hand side of \eqref{eq:uni3}, invoking \eqref{eq:uni1}, we have
\begin{equation}
|Re(\kappa+i\zeta)(g,w)_h|\leqslant \frac{3}{2}\sqrt{\kappa^2+\zeta^2}(1+C_u)^2\|w\|^2_h.
\end{equation}
Substituting above inequality into \eqref{eq:uni3} yields
\begin{equation}
\|w\|^2_h\leqslant \t \frac{3\sqrt{\kappa^2+\zeta^2}(1+C_u)^2+|\gamma|}{2}\|w\|^2_h.
\end{equation}
For $\t<\frac{2}{3\sqrt{\kappa^2+\zeta^2}(1+C_u)^2+|\gamma|}$, we obtain $\|w\|_h=0$, which implies
\begin{equation}
w_j=0, \quad j\in\mathcal{T}_M.
\end{equation}
Thus, the proof is complete.
$\Box$
\section{ Numerical experiments}
%
In this section, we report some numerical results to confirm the
theoretical accuracy and  efficiency of scheme \eqref{eq:CCNS}-\eqref{eq:DBC}.
\subsection{Iterative algorithm}
%
Before embarking on our numerical experiments, developing
an efficient iterative algorithm is of the essence to compute the solution
of the system of nonlinear equation arising at a given time level of the scheme \eqref{eq:CCNS}-\eqref{eq:DBC}.

To this end, rewrite \eqref{eq:CCNS} as
\begin{equation}
u^{n+\frac{1}{2}}_j=u^n_j-\frac{\tau}{2}\Big[(\upsilon+i\eta)\Delta^{\a}_{h}u^{n+\frac{1}{2}}_j+(\kappa+i\zeta)|u^{n+\frac{1}{2}}_j|^2u^{n+\frac{1}{2}}_j
-\gamma u^{n+\frac{1}{2}}_j\Big].
\end{equation}
We introduce the variable $z_j=u^{n+\frac{1}{2}}_j$
and then obtain
\begin{equation}
z_j=u^n_j-\frac{\tau}{2}\Big[(\upsilon+i\eta)\Delta^{\a}_{h}z_j+(\kappa+i\zeta)|z_j|^2z_j
-\gamma z_j\Big],
\end{equation}
and
\begin{equation}\label{eq:uz}
u^{n+1}_j=2z_j-u^n_j.
\end{equation}
Then we propose the following iterative algorithm
 \begin{equation}\label{eq:ITSC}
 \begin{split}
  &z^{(s+1)}_j=u^n_j-\frac{\tau}{2}\Big[(\upsilon+i\eta)\Delta^{\a}_{h}z^{(s+1)}_j+(\kappa+i\zeta)|z^{(s)}_j|^2z^{(s)}_j-\gamma z^{(s+1)}_j\Big],\\
  &\ \ \ \ \ \ \ \ \ \ \ \ \ \ \ \ \ \ \ \ \ \ \ \ \ \ \ \ \ \ \ \ \ \ \ \ \ \ \ \ \ \ j\in\mathcal{T}_M, \ 0\leqslant n \leqslant N-1, \ s=0,1,2,\ldots,\\
  &z^{(s)}_j=0, \quad j\in\mathbb{Z}\backslash\mathcal{T}_M, \quad 0\leqslant n \leqslant N,\ s=0,1,2,\ldots.
 \end{split}
 \end{equation}
The initial iteration is selected as, for $n\geqslant 1$,
\[
z^{(0)}_j=\frac{3}{2}u^n_j-\frac{1}{2}u^{n-1}_j,
\]
and for $n=0$,
\[
z^{(0)}_j=u^0_j-\frac{\tau}{2}\big[(\upsilon+i\eta)\Delta^{\a}_{h}u^0_j+(\kappa+i\zeta)|u^0_j|^2u^0_j-\gamma u^0_j\big].
\]

The system is indeed linearized at each iteration, and we solve an
inner problem $Az^{(s+1)}=b$ to get $z^{(s+1)}_j$. Then $z_j$ is numerically reached once $z^{(s+1)}_j$
converges and $u^{n+1}_j$ is obtained by \eqref{eq:uz}. It is noted that the coefficient matrix $A$ is
independent of the time level, and this feature dramatically
benefits the numerical implementation.
\subsection{Numerical tests}
\noindent\textbf{Example 5.2.1}\ \
We firstly testify the numerical accuracy of the scheme with
\[
\upsilon=0.3, \eta=\frac{1}{2}, \kappa=-\frac{\upsilon(3\sqrt{1+4\upsilon^2}-1)}{2(2+9\upsilon^2)}, \zeta=-1, \gamma=0.
\]
For $\a=2$, the exact solution is explicitly given by \cite{AkhmedievNN:1996:PRE}
\begin{equation}\label{eq:exs}
u(x,t)=a(x)\exp{(id\ln[a(x)]-i\omega t)},
\end{equation}
where
\[
a(x)=F\text{sech}{(x)}, F=\sqrt{\frac{d\sqrt{1+4\upsilon^2}}{-2\kappa}}, d=\frac{\sqrt{1+4\upsilon^2}-1}{2\upsilon}, \omega=-\frac{d(1+4\upsilon^2)}{2\upsilon}.
\]
The computational interval is chosen as $[a,b]=[-16,16]$ and the initial value is taken as $u(x,0)$ in \eqref{eq:exs}. Choose the iteration tolerance as $10^{-14}$ for the iterative
algorithm \eqref{eq:ITSC}.
For $1<\a<2$, the exact solution can not be explicitly given and thus, the numerical ``exact'' solution $u$ is computed using the proposed scheme with a very fine mesh size $h=0.0125$
and time step $\t=0.0001$. Let $u_h$ be the numerical solution. We
measure the error $e(\t,h)=u-u_h$ at time $T=1$ with the $l^2$ norm
and the $l^\infty$ norm. The corresponding convergence
orders are calculated by
\begin{equation*}
\begin{split}
&\text{order}1=\log_2(\|e(h,\t)\|_h/\|e(\t/2,h/2)\|_h),\\
&\text{order}2=\log_2(\|e(h,\t)\|_{l^\infty_h}/\|e(\t/2,h/2)\|_{l^\infty_h}).
\end{split}
\end{equation*}
Tables \ref{tab:order1}-\ref{tab:order2} list the errors and corresponding orders of the numerical
scheme for $\a=2$ and $1<\a<2$, respectively.
The data confirm the theoretical accuracy of the difference
scheme \eqref{eq:CCNS}-\eqref{eq:DBC} in Theorem \ref{th:Converg}.

\begin{table}[!h]
\tabcolsep 0pt \caption{The $l^2_h$ and $l^{\infty}_h$ errors and their convergence orders for $\a=2$.}
\vspace*{-12pt}
\begin{center}
\def\temptablewidth{0.8\textwidth}
{\rule{\temptablewidth}{1pt}}
\begin{tabular*}{\temptablewidth}{@{\extracolsep{\fill}}llllll}
 $\tau$     & $h$         &  $\|e\|_h$      & order1   & $\|e\|_{l^{\infty}_h}$  & order2    \\   \hline\noalign{\smallskip}
 0.02       & 0.2         &  5.5462e-003    &  $-$     & 5.5486e-003       & $-$       \\
 0.01       & 0.1         &  1.3766e-003    &  2.0104  & 1.3691e-003       & 2.0190  \\
 0.005      & 0.05        &  3.4353e-004    &  2.0026  & 3.4117e-004       & 2.0046  \\
 0.0025     & 0.025       &  8.5845e-005    &  2.0006  & 8.5225e-005       & 2.0011  \\
 0.00125    & 0.0125      &  2.1460e-005    &  2.0001  & 2.1302e-005       & 2.0003  \\
\end{tabular*}
{\rule{\temptablewidth}{1pt}}\label{tab:order1}
\end{center}
\end{table}
\begin{table}[!h]
\tabcolsep 0pt \caption{The $l^2_h$ and $l^{\infty}_h$ errors and their convergence orders for $1<\a<2$.}
\vspace*{-12pt}
\begin{center}
\def\temptablewidth{0.9\textwidth}
{\rule{\temptablewidth}{1pt}}
\begin{tabular*}{\temptablewidth}{@{\extracolsep{\fill}}lllllll}
 $\a$     & $h$        & $\t$     & $\|e\|_h$    & order1   & $\|e\|_{l^{\infty}_h}$  & order2 \\ \hline\noalign{\smallskip}
 {}       & 0.2        & 0.02     & 1.2966e-002  & $-$      & 1.8415e-002        & $-$      \\
 1.3      & 0.1        & 0.01     & 3.1803e-003  & 2.0275   & 4.4581e-003        & 2.0464 \\
 {}       & 0.05       & 0.005    & 7.5488e-004  & 2.0749   & 1.0528e-003        & 2.0822 \\ \hline\noalign{\smallskip}
 {}       & 0.2        & 0.02     & 1.0519e-002  & $-$      & 1.3001e-002        & $-$      \\
 1.6      & 0.1        & 0.01     & 2.5499e-003  & 2.0444   & 3.0928e-003        & 2.0717 \\
 {}       & 0.05       & 0.005    & 6.0393e-004  & 2.0780   & 7.2919e-004        & 2.0846 \\ \hline\noalign{\smallskip}
 {}       & 0.2        & 0.02     & 6.7430e-003  & $-$      & 7.0782e-003        & $-$      \\
 1.9      & 0.1        & 0.01     & 1.6458e-003  & 2.0346   & 1.7127e-003        & 2.0471 \\
 {}       & 0.05       & 0.005    & 3.9052e-004  & 2.0753   & 4.0552e-004        & 2.0784 \\
\end{tabular*}
{\rule{\temptablewidth}{1pt}}\label{tab:order2}
\end{center}
\end{table}
\noindent\textbf{Example 5.2.2}\ \
In this tests, the initial value is chosen as
\[
u(x,0)=\exp(-2x^2).
\]
Henceforth, we always take $[a,b]=[-10,10]$ and $\t=h=0.05$.

\begin{figure}[!ht]
\centering \subfigure{\includegraphics[width=8cm]{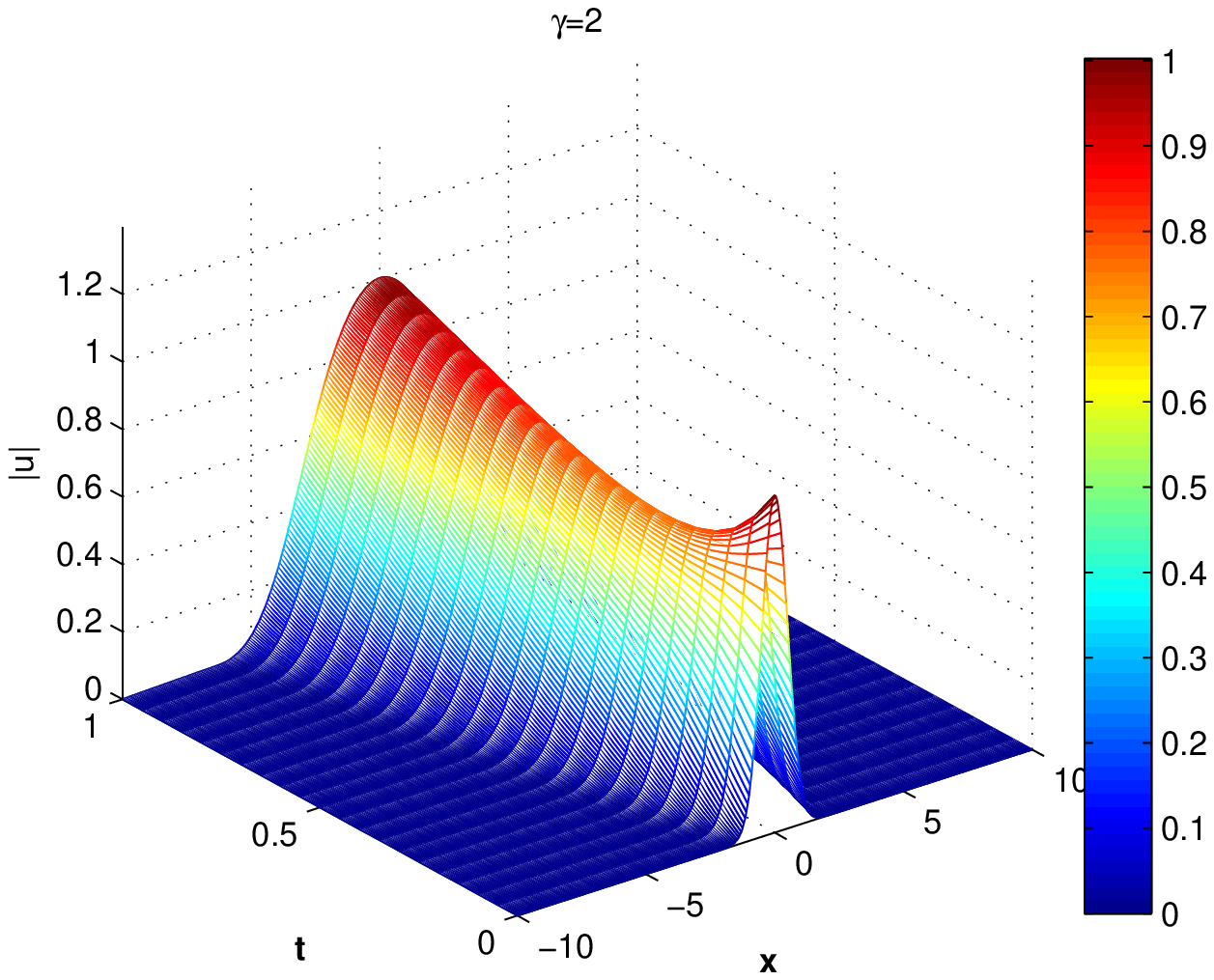}}
\subfigure{\includegraphics[width=8cm]{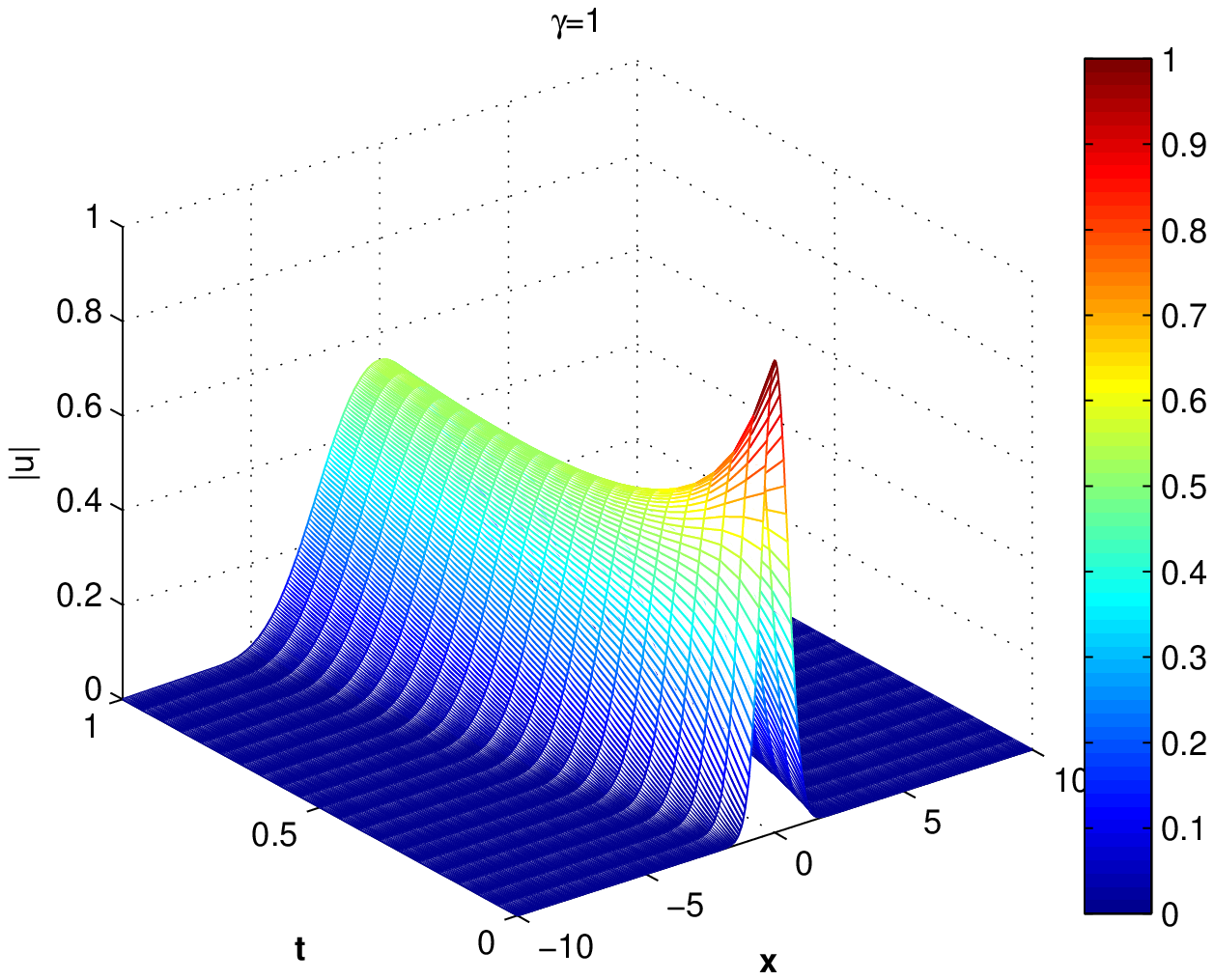}} \caption{The evolution of $|u|$ for $\gamma=2$ (left) and $\gamma=1$ (right).} \label{fig:g1}
\end{figure}

\begin{figure}[!ht]
\centering \subfigure{\includegraphics[width=8cm]{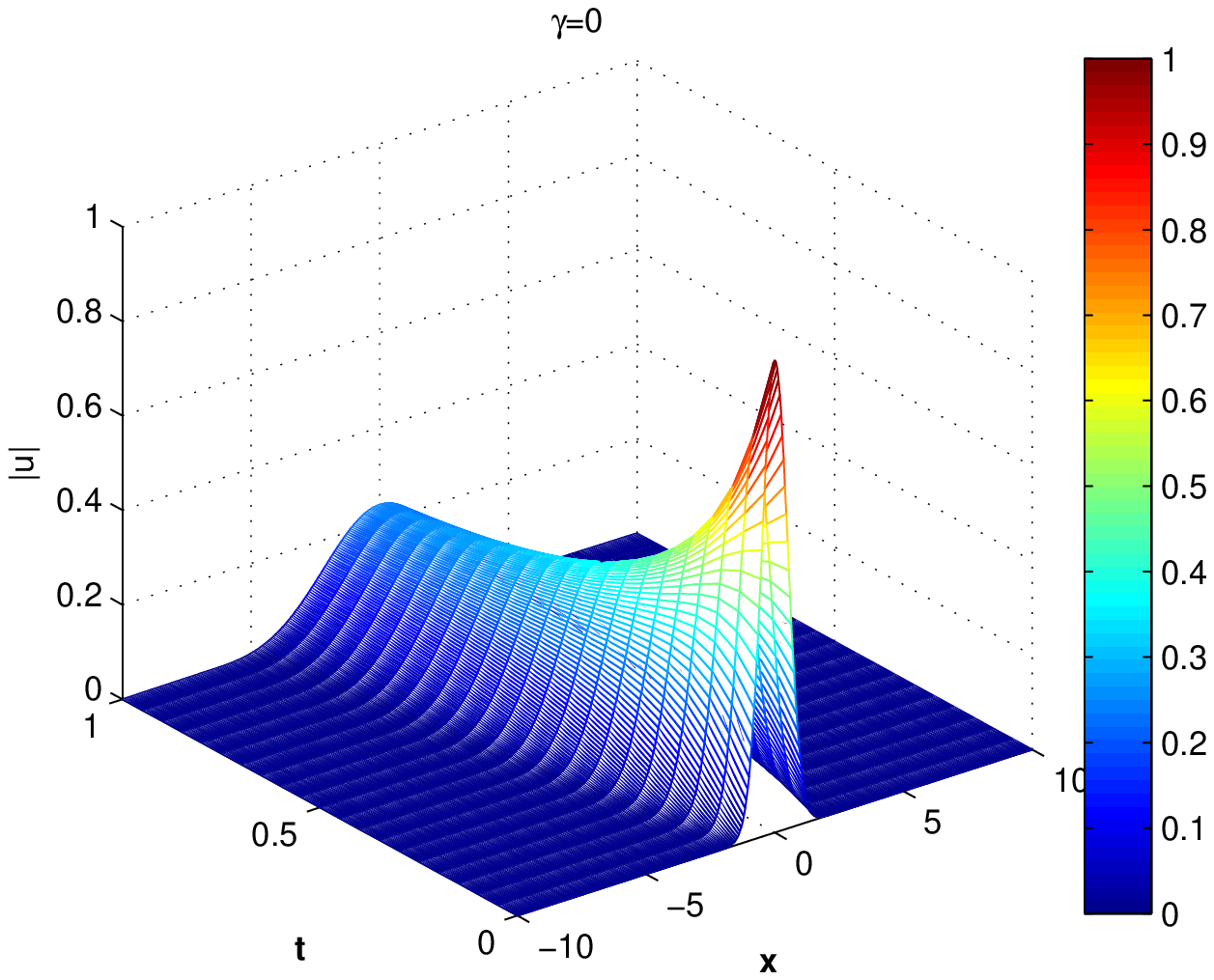}}
\subfigure{\includegraphics[width=8cm]{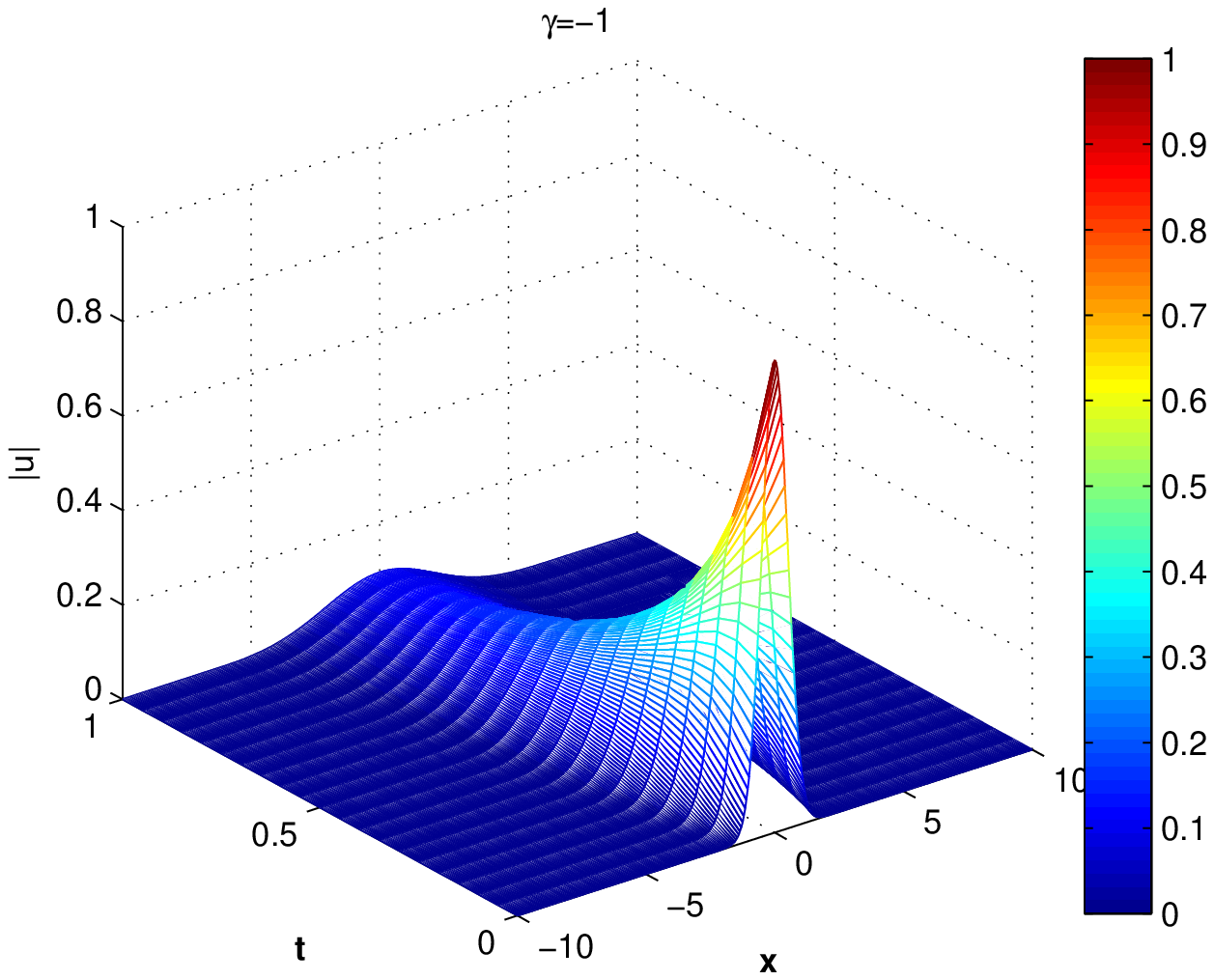}} \caption{The evolution of $|u|$ for $\gamma=0$ (left) and $\gamma=-1$ (right).} \label{fig:g2}
\end{figure}

\begin{figure}[!ht]
\centering
\includegraphics[width=9cm]{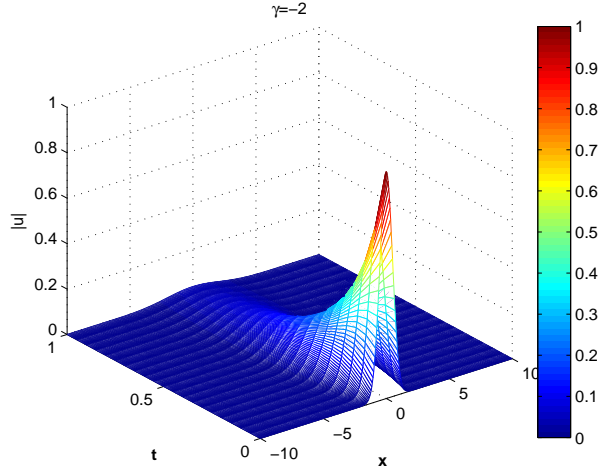}
\caption{The evolution of $|u|$ for $\gamma=-2$.} \label{fig:g3}
\end{figure}

\begin{figure}[!ht]
\centering \subfigure{\includegraphics[width=8cm]{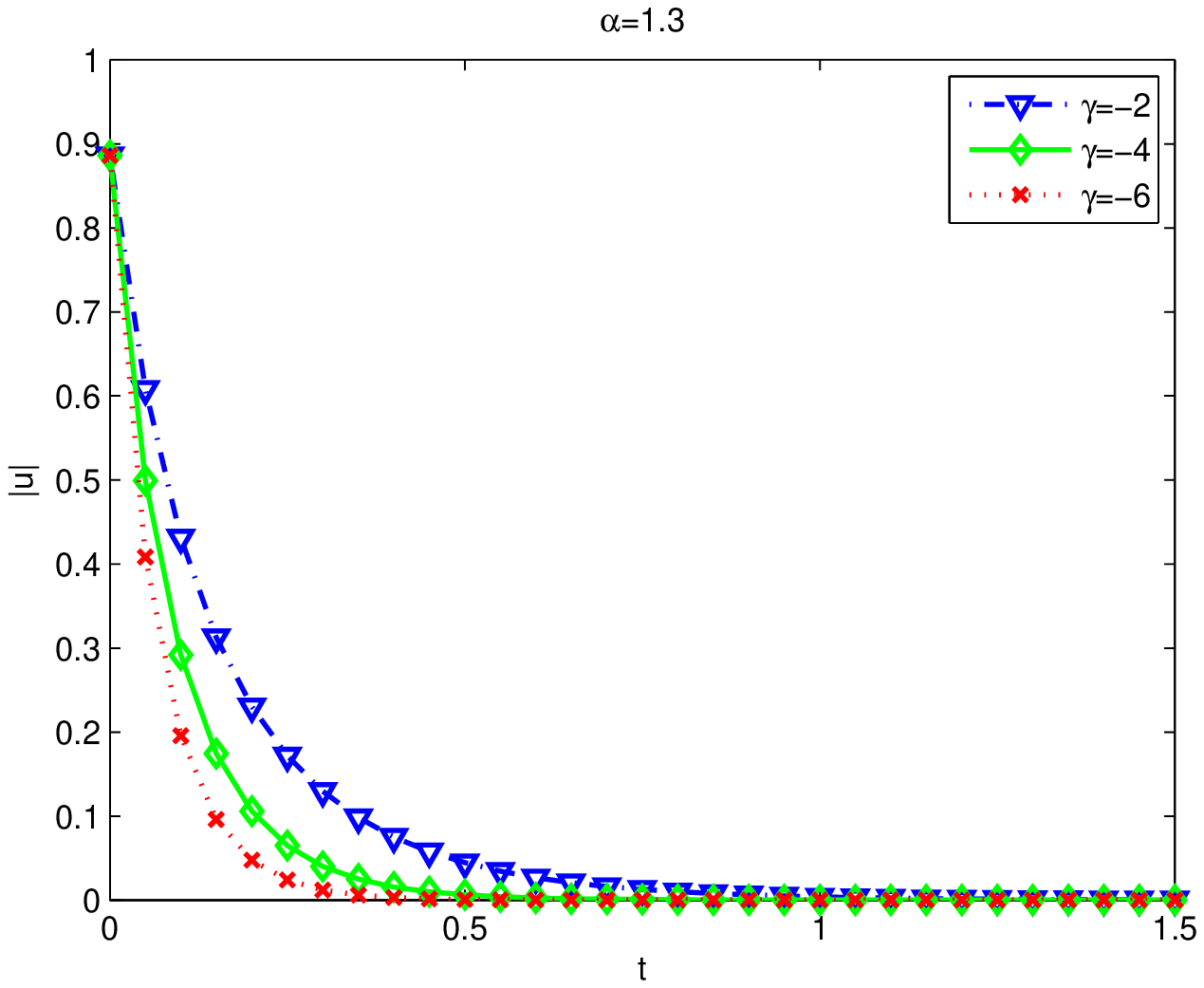}}
\subfigure{\includegraphics[width=8cm]{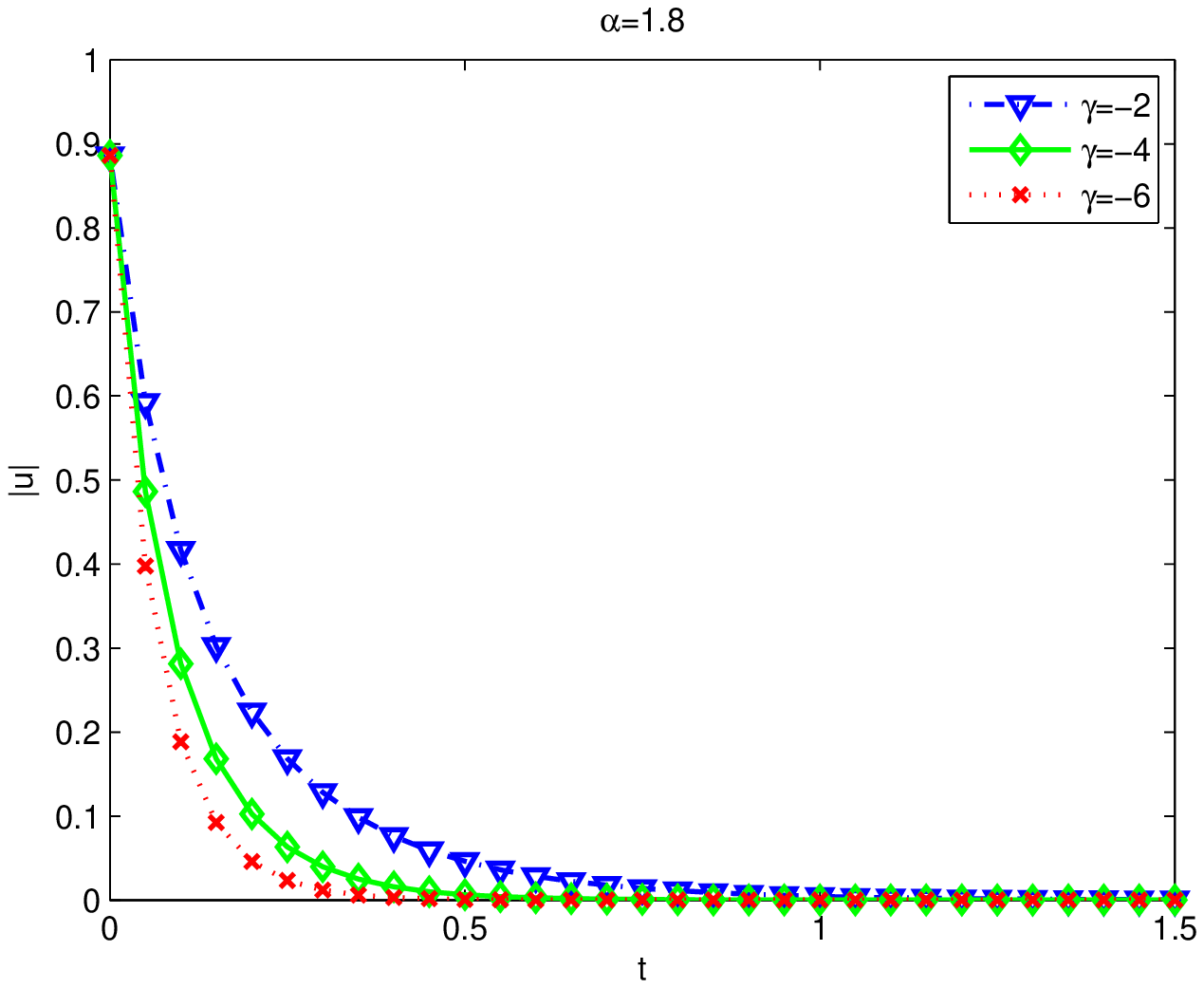}} \caption{The discrete norm $\|u\|^2_h$ at for $\a=1.3$ (left) and $\a=1.8$ (right) with $\gamma=-2,-4,-6$.} \label{fig:g4}
\end{figure}

Firstly, the evolution of the numerical solution is depicted.
We pay particular attention to the influence of parameter $\gamma$ on the
evolution of wave-shape in the fractional case. To this end, chosen
$\upsilon=1, \eta=1, \kappa=1, \zeta=2$, we take different values of
$\gamma$, i.e., $\gamma=2,1,0,-1,-2$ with $\a=1.8$ and depict the
evolution of $|u|$ in Figures \ref{fig:g1}-\ref{fig:g3}. It is
observed that, as in the classical case, the parameter $\gamma$
dramatically affects the wave-shape in the fractional case. In
addition, the solution decays rapidly with time evolution especially
for $\gamma<0$. For more intensive study, then in Figure
\ref{fig:g4}, we further depict the evolution of $\|u\|^2_h$ with
$\alpha=1.3, 1.8$. Recall that, in the classical case ($\a=2$), the
discrete norm $\|u\|^2_h$ decays to zero for $\gamma<0$ and when
$\gamma$ is smaller, the decay is faster \cite{XuQ:2011:NMPDE}. In
our fractional case, we observe similar phenomena and the fractional
order $\a$ affects the decay rate very slightly.

\begin{figure}[!ht]
\centering
\includegraphics[width=9cm]{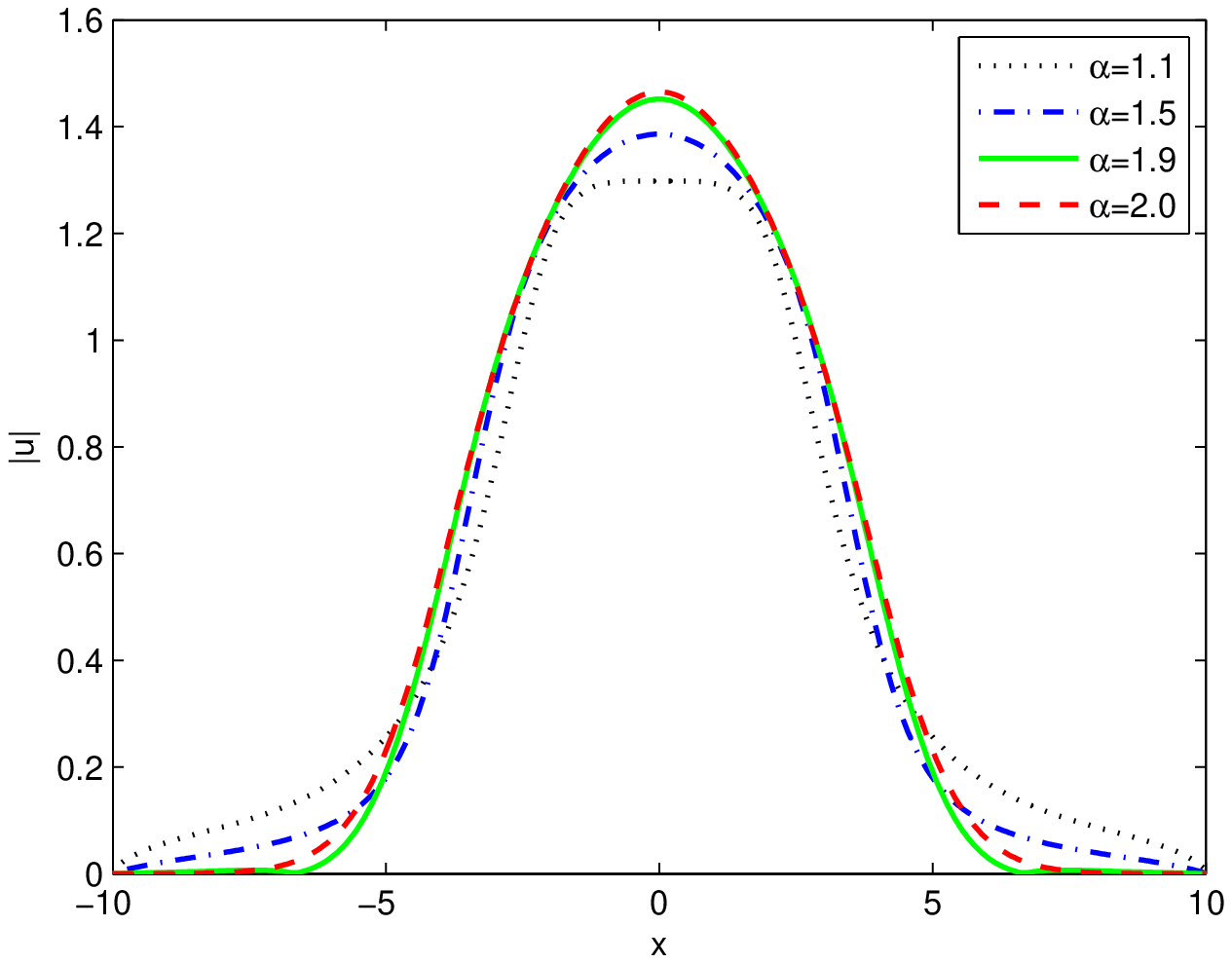}
\caption{The profile of $|u|$ at $t=1$ with $\a=1.1, 1.5, 1.9, 2.0$.} \label{fig:a1}
\end{figure}

Secondly, we numerically study the impact on the dissipative mechanism of the fractional
Laplacian. Choose $\upsilon=1, \eta=1, \kappa=1, \zeta=2, \gamma=3$.
The solutions at $t=1$ with different values of $\a$ are depicted in Figure \ref{fig:a1}.
It can be seen that the wave-shape changes with fractional parameter $\a$. This phenomenon is greatly different
from that in the classical case and essentially, features the nonlocal character
of the fractional Laplacian.

\begin{figure}[!ht]
\centering
\subfigure{\includegraphics[width=8cm]{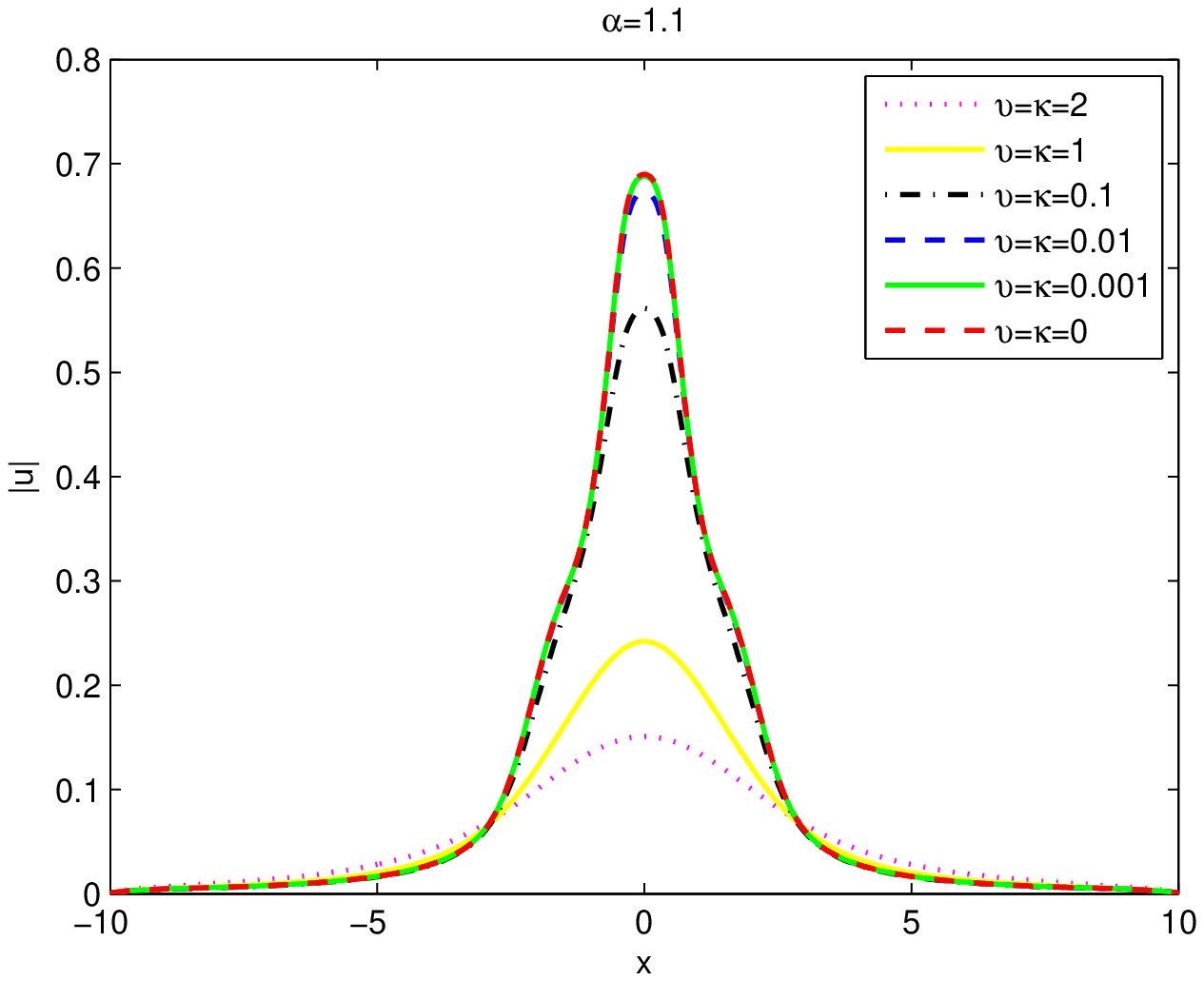}}
\subfigure{\includegraphics[width=8cm]{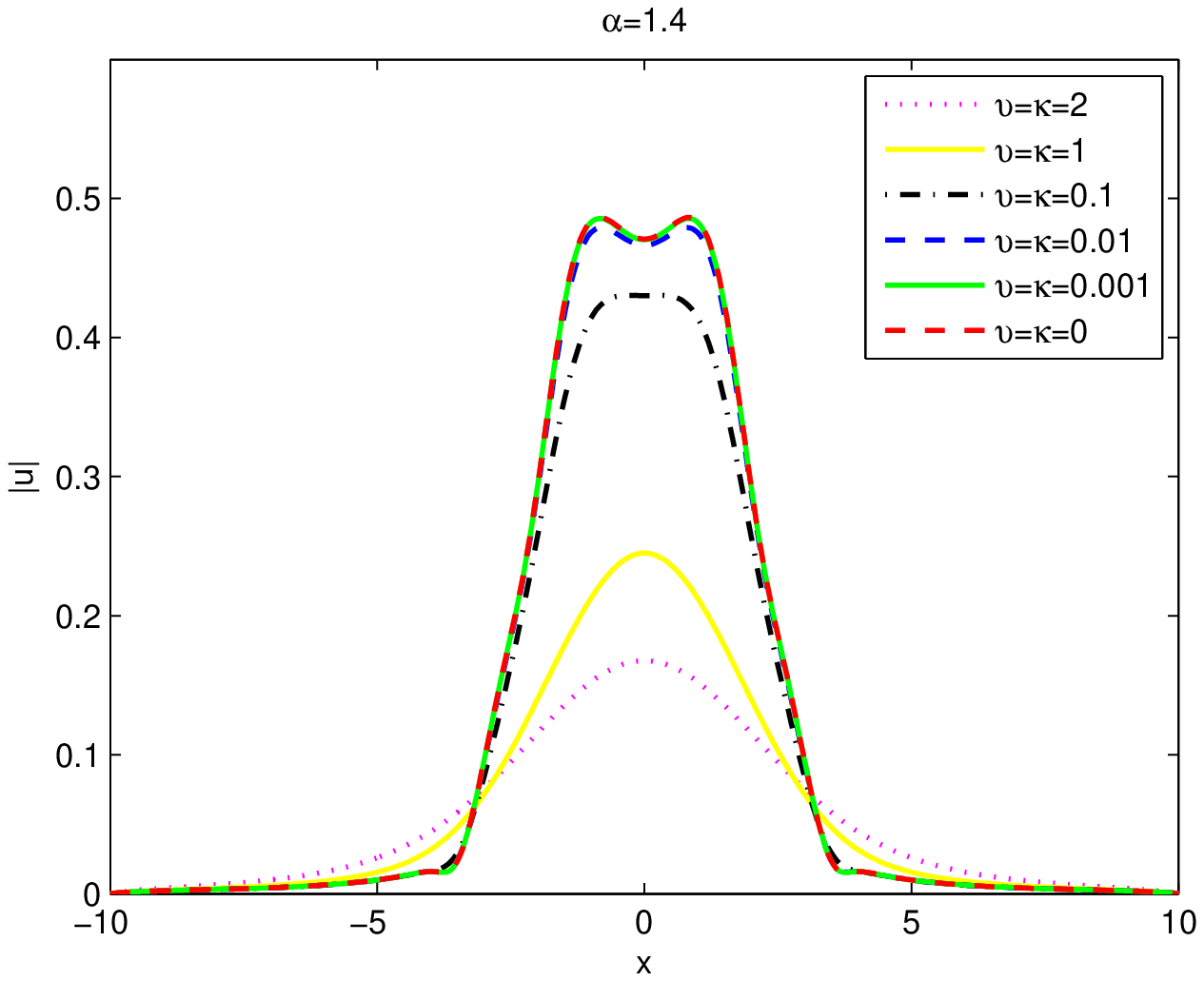}}
\subfigure{\includegraphics[width=8cm]{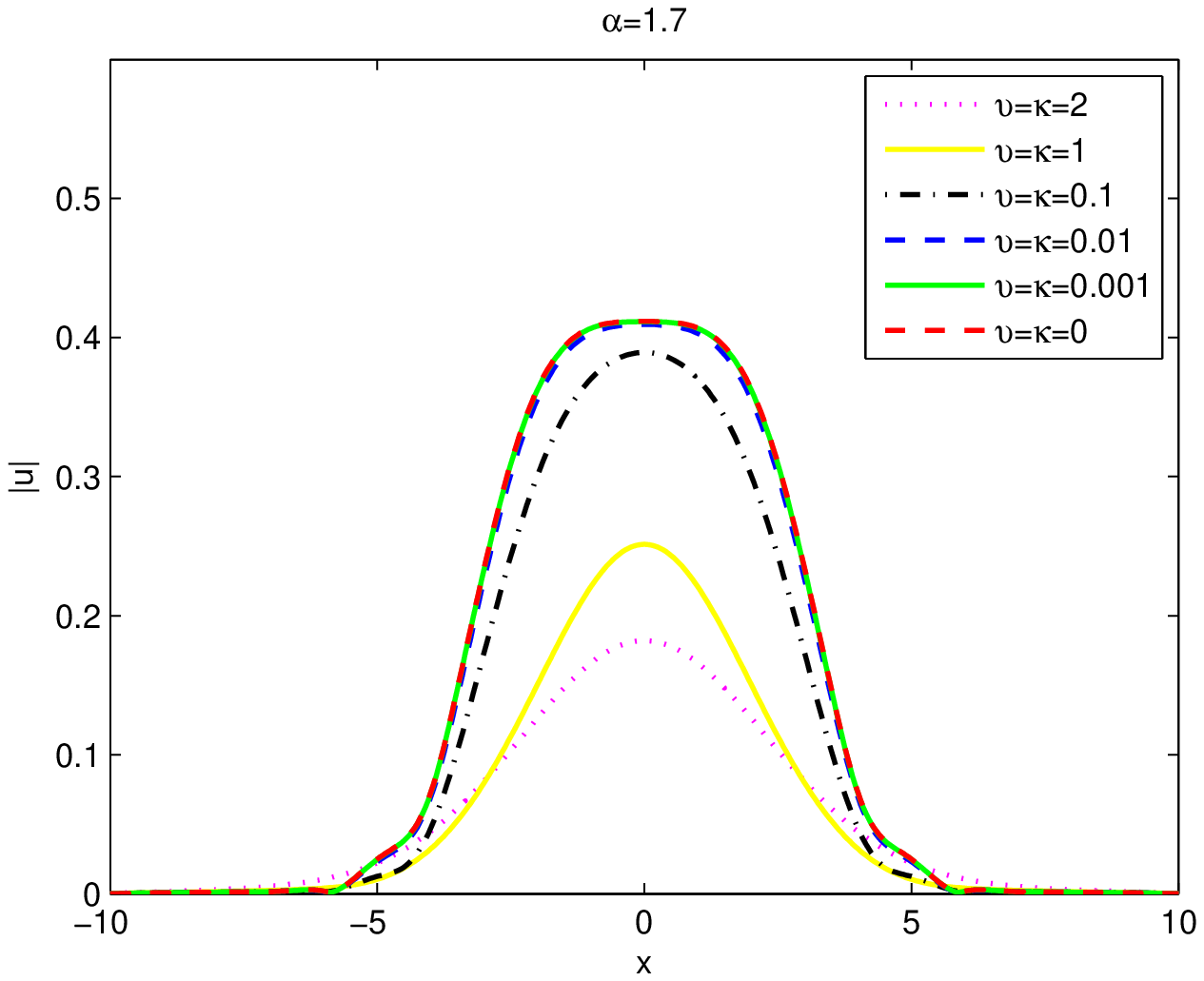}}
\subfigure{\includegraphics[width=8cm]{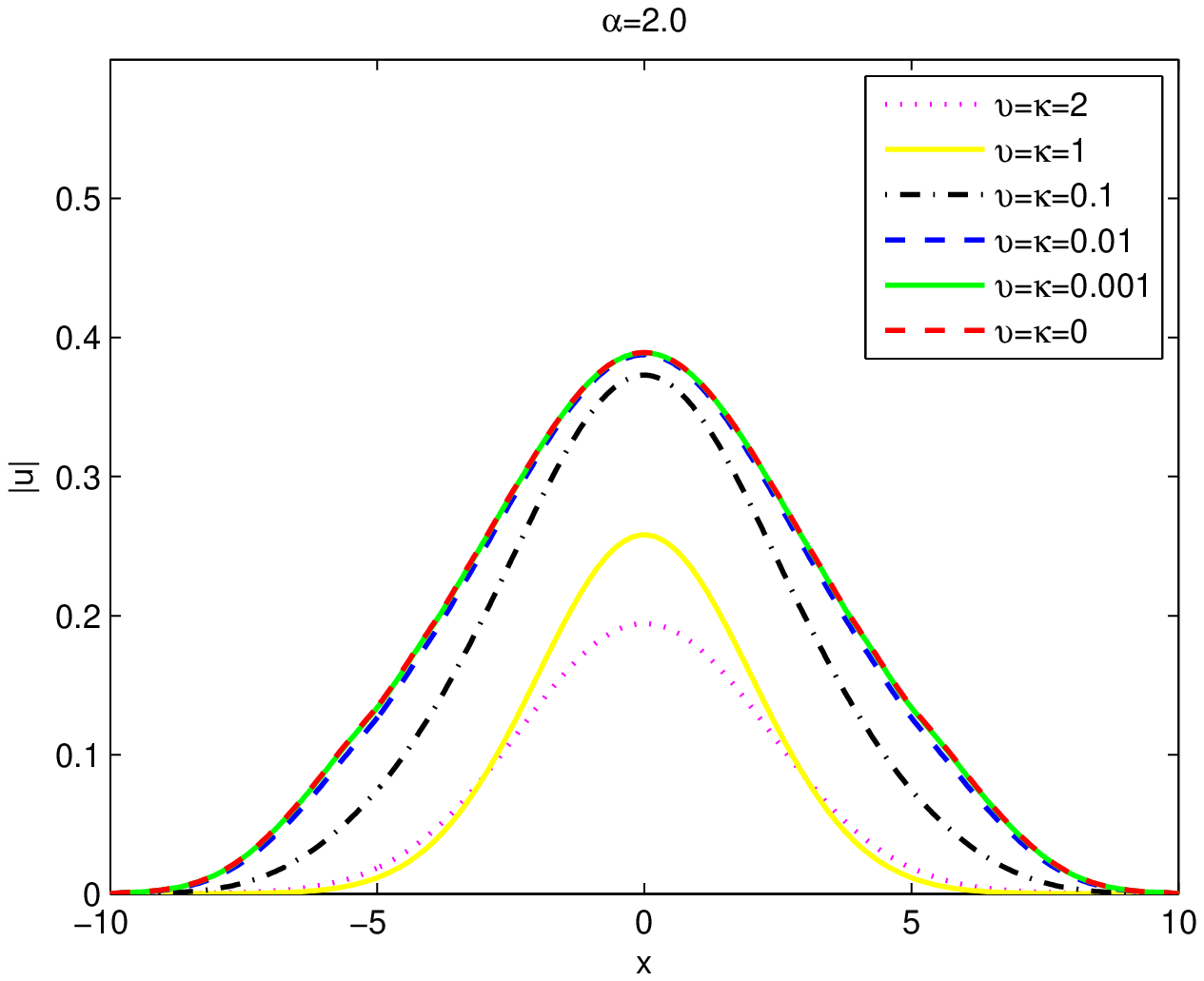}}
\caption{The profile of $|u|$ at $t=1$ with diminishing $\upsilon$ and $\kappa$ for $\a=1.1, 1.4, 1.7, 2.0$.} \label{fig:v4}
\end{figure}

Finally, we simulate the inviscid limit behavior of the solution. The
authors in \cite{GuoB:2013:FCAA} have shown that the solution of the
FGLE converges to the solution of the FSE (i.e.,
$\upsilon=0,\kappa=0$) when $\upsilon\rightarrow 0,\kappa\rightarrow
0$. For the numerical simulation of the FSE, see, e.g.,
\cite{WangP:2015:JCP,WangP:2015:NA,WangD:2013:JCP,WangD:2014:JCP}.
Here we numerically testify this matter. For this purpose, we set
$\eta=1, \zeta=-2, \gamma=0$ and choose diminishing $\upsilon$ and
$\kappa$. From Figure \ref{fig:v4}, it is observed that the solution
asymptotically approachs to the solution of the FSE. This
observation confirms the theoretical results in
\cite{GuoB:2013:FCAA}.
%
\section{ Conclusions}
%
In this paper, we proposed and analyzed a finite difference scheme
for solving the nonlinear complex fractional Ginzburg-Landau
equation where the fractional Laplacian was approximated by the
weighted and shifted Gr\"{u}nwald difference operator. We obtained
the unconditional optimal convergence rate at the order of
$O(\t^2+h^2)$ in the $l^2_h$ norm with the time step $\t$ and mesh
size $h$, without any mesh ratio constraints. In the proof for the
scheme, building on the careful analysis of the difference operator,
we established some useful inequalities with respect to the
fractional Sobolev norm and the a priori bound of the numerical
solution. Both theoretical analysis and numerical tests show that
the scheme is efficient for the numerical solution of the nonlinear
fractional Ginzburg-Landau equation.
%

%

\section*{Acknowledgments}
%
This work was supported by National Natural Science Foundation of
China (Nos. 91130003 and 11371157) and the Graduates' Innovation
Fund of Huazhong University of Science \& Technology (No.
2015650011). The authors wish to thank the
anonymous referees for their valuable comments and suggestions which
lead to an improvement of this paper.

\medskip

\medskip

\noindent\textbf{Appendix A: Hausdorff-Young inequality}\\

\noindent\textbf{Lemma A.1.} If $1\leqslant q \leqslant 2$, $\frac{1}{q}+\frac{1}{p}=1$, then
\begin{equation}
\big(h\sum_{j\in\mathbb{Z}}|u_j|^p\big)^{\frac{1}{p}}\leqslant C \big(\int^{\pi/h}_{-\pi/h}|\widehat{u}(k)|^q\text{d}k\big)^{\frac{1}{q}},
\end{equation}
where $\frac{1}{\sqrt{2\pi}}\leqslant C\leqslant 1$.\\

\noindent\textbf{Proof.} The proof is very similar to that of the
same inequality for Fourier series (see Corollary 2.4 in \cite[page
57]{SteinE:2011:PUR}) and hence, we just give below a sketch of it.
From the definition of inversion formula, we have
\[
\sup_{j\in\mathbb{Z}}|u_j|\leqslant \frac{1}{\sqrt{2\pi}} \int^{\pi/h}_{-\pi/h}|\widehat{u}(k)|\text{d}k.
\]
The Parseval's indentity gives
\[
h\sum_{j\in\mathbb{Z}}|u_j|^2=\int^{\pi/h}_{-\pi/h}|\widehat{u}(k)|^2\text{d}k.
\]
Then invoking the Riesz interpolation theorem (see Theorem 2.1 in
\cite[page 52]{SteinE:2011:PUR}), we obtain the conclusion. $\Box$

\medskip

\medskip

%
\noindent\textbf{Appendix B: The proof of Lemma \ref{Lem:L1}}\\

\noindent\textbf{Lemma \ref{Lem:L1}.}
For $1<\a\leqslant 2$, let $h(\a,\omega)$ be the function defined by
\[
h(\a,\omega)=\lambda_1\cos\big(\frac{\a}{2}(\omega-\pi)-\omega\big)
+\lambda_0\cos\big(\frac{\a}{2}(\omega-\pi)\big)+\lambda_{-1}\cos\big(\frac{\a}{2}(\omega-\pi)+\omega\big),
\]
where $\omega\in[0,\pi]$ and $\lambda_1, \lambda_0, \lambda_{-1}$ are defined in \eqref{eq:c1}.
Then $h(\a,\omega)$ does not decrease with respect to $\omega$.\\

\noindent\textbf{Proof.} It is easy to see that $h(\a,\omega)\equiv-1$ for $\a=2$.

For $1<\a< 2$, rearranging $h(\a,\omega)$ gives
\begin{equation*}
\begin{split}
h(\a,\omega)&=\frac{\a^2+3\a+2}{12}\cos\big(\frac{\a}{2}(\omega-\pi)-\omega\big)+\frac{4-\a^2}{6}\cos\big(\frac{\a}{2}(\omega-\pi)\big)\\
&\ \ \ \ +\frac{\a^2-3\a+2}{12}\cos\big(\frac{\a}{2}(\omega-\pi)+\omega\big)\\
&=\frac{\a^2+2}{12}\Big(\cos\big(\frac{\a}{2}(\omega-\pi)-\omega\big)+\cos\big(\frac{\a}{2}(\omega-\pi)+\omega\big)\Big)\\
&\ \ \ \ +\frac{3\a}{12}\Big(\cos\big(\frac{\a}{2}(\omega-\pi)-\omega\big)-\cos\big(\frac{\a}{2}(\omega-\pi)+\omega\big)\Big)+\frac{4-\a^2}{6}\cos\big(\frac{\a}{2}(\omega-\pi)\big)\\
&=\frac{\a^2+2}{6}\cos\big(\frac{\a}{2}(\omega-\pi)\big)\cos(\omega)+\frac{\a}{2}\sin\big(\frac{\a}{2}(\omega-\pi)\big)\sin(\omega)+\frac{4-\a^2}{6}\cos\big(\frac{\a}{2}(\omega-\pi)\big).
\end{split}
\end{equation*}
Taking the derivative of $h(\a,\omega)$ with respect to $\omega$, we obtain
\begin{equation*}
\begin{split}
h'(\a,\omega)&=\frac{\a^2+2}{6}\Big(-\frac{\a}{2}\sin\big(\frac{\a}{2}(\omega-\pi)\big)\cos(\omega)-\cos\big(\frac{\a}{2}(\omega-\pi)\big)\sin(\omega)\Big)\\
&\ \ \ \ +\frac{\a}{2}\Big(\frac{\a}{2}\cos\big(\frac{\a}{2}(\omega-\pi)\big)\sin(\omega)+\sin\big(\frac{\a}{2}(\omega-\pi)\big)\cos(k)\Big)
-\frac{\a(4-\a^2)}{12}\sin\big(\frac{\a}{2}(\omega-\pi)\big)\\
&=\frac{\a(4-\a^2)}{12}\sin\big(\frac{\a}{2}(\omega-\pi)\big)\big(\cos(\omega)-1\big)-\frac{4-\a^2}{12}\cos\big(\frac{\a}{2}(\omega-\pi)\big)\sin(\omega)\\
&=-\frac{\a(4-\a^2)}{6}\sin\big(\frac{\a}{2}(\omega-\pi)\big)\sin^2\big(\frac{\omega}{2}\big)
-\frac{4-\a^2}{6}\cos\big(\frac{\a}{2}(\omega-\pi)\big)\sin\big(\frac{\omega}{2}\big)\cos\big(\frac{\omega}{2}\big)\\
&=-\frac{4-\a^2}{6}\sin\big(\frac{\omega}{2}\big)g(\a,\omega),
\end{split}
\end{equation*}
where
\[
g(\a,\omega):=\a\sin\big(\frac{\a}{2}(\omega-\pi)\big)\sin\big(\frac{\omega}{2}\big)
+\cos\big(\frac{\a}{2}(\omega-\pi)\big)\cos\big(\frac{\omega}{2}\big).
\]
In order to show that $h(\a,\omega)$ non-decreases with respect to $\omega$, it is sufficient to show that $h'(\a,\omega)\geqslant 0$, or $g(\a,\omega)\leqslant 0$ for $\omega\in[0,\pi]$. In fact, taking the derivative of $g(\a,\omega)$ with respect to $\omega$ yields
\begin{equation}
g'(\a,\omega)=\frac{\a^2-1}{2}\cos\big(\frac{\a}{2}(\omega-\pi)\big)\sin\big(\frac{\omega}{2}\big),
\end{equation}
which implies that the extreme point of $g(\a,\omega)$ is $\omega=\frac{\a-1}{\a}\pi$ for $\omega\in(0,\pi)$. By a simple analysis, it is shown that $g(\a,\omega)$ reach its maximum value at $\omega=\pi$ and
reach its minimum value at $\omega=\frac{\a-1}{\a}\pi$, which gives
\[
-\a\sin\big(\frac{\a-1}{2\a}\pi\big)=g(\a,\frac{\a-1}{\a}\pi)\leqslant g(\a,\omega)\leqslant g(\a,\pi)=0.
\]
Thus, the proof is complete.
$\Box$

\bibliographystyle{model1-num-names}
\bibliography{Schrodinger1}



\end{document}